\documentclass[11pt, reqno]{amsart}
\usepackage{amsmath}
\usepackage{amssymb}
\usepackage{amstext}
\usepackage{mathrsfs}
\usepackage{a4wide}
\usepackage{graphicx}
\usepackage{verbatim}
\usepackage{bm}

\usepackage{color}
\usepackage{cases}
\usepackage{hyperref}
\hypersetup{colorlinks=true, linkcolor=blue, anchorcolor=blue, citecolor=blue }

\newtheorem{theorem}{Theorem}[section]
\newtheorem{proposition}[theorem]{Proposition}
\newtheorem{corollary}[theorem]{Corollary}
\newtheorem{lemma}[theorem]{Lemma}
\theoremstyle{remark}
\newtheorem{remark}{Remark}[section]

\theoremstyle{definition}
\numberwithin{equation}{section}

\newcommand{\R}{\mathbb{R}}

\allowdisplaybreaks

\title[Norbury vortex rings]{Global continuation of Norbury vortex rings}
\author{Guolin Qin}
\address{State Key Laboratory of Mathematical Sciences, Academy of Mathematics and Systems Science, Chinese Academy of Sciences, Beijing 100190, P.R. China}
\email{qinguolin18@mails.ucas.ac.cn}

\author{Weicheng Zhan}
\address{School of Mathematical Sciences, Xiamen University, Xiamen, Fujian 361005, P.R. China}
\email{zhanweicheng@amss.ac.cn}

\begin{document}

\begin{abstract}
We study Norbury's family of vortex rings~\cite{Nor73} for the axisymmetric Euler equations without swirl and prove that the continuum from Hill's spherical vortex constructed by Amick and Turner~\cite{AT88} reaches every nonnegative flux. We first bound the energy and the support of all centered solutions with bounded flux. Using this bound and the result of Amick and Fraenkel~\cite{AF88}, we prove uniqueness among centered finite-energy solutions for sufficiently small flux. As the flux tends to infinity, the core moves away from the axis and its rescaled cross-section converges to a disk. We further show that the centered solution at sufficiently large flux is unique and depends continuously on the flux.
\end{abstract}

\maketitle

\section{Introduction}\label{sec1}

We consider steadily translating vortex rings for the three-dimensional incompressible Euler equations
\begin{equation*}
	\partial_tu+(u\cdot\nabla)u+\nabla p=0,\qquad
	\nabla\cdot u=0\qquad\text{in }\R^3,
\end{equation*}
where $u$ is the velocity and $p$ is the pressure. The mathematical study of vortex rings goes back to Helmholtz \cite{Hel} and to the work of Kelvin and Hicks \cite{Tho,Hicks85}. We refer to \cite{Lam,Saf,Akh,MGT} for accounts of the classical theory. For experimental studies of vortex rings and their role in propulsion, see \cite{Max,LT,Sa92}.

\begin{figure}[t]
	\centering
	\includegraphics[width=\linewidth]{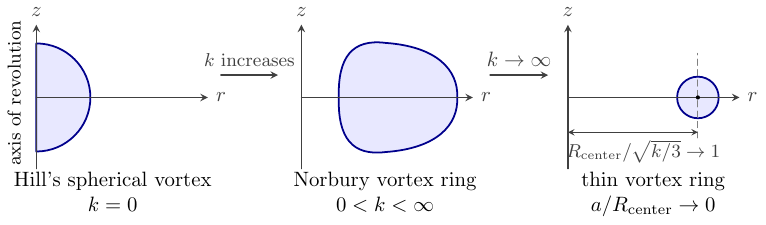}
	\caption{Meridional cross-sections along the Norbury family. Here $r$ and $z$ are the radial and axial coordinates, and $k$ is the flux. The radial coordinate of the center of the section is $R_{\mathrm{center}}$, and $a$ is the radius of a disk with the same area as the section. The middle drawing is schematic.}
	\label{fig1-1}
\end{figure}

For axisymmetric flows without swirl we take the ratio of the azimuthal vorticity to the distance from the symmetry axis to be a positive constant in the core and zero outside it. This choice includes Hill's spherical vortex \cite{Hil}, the toroidal rings close to it constructed by Norbury \cite{Nor72}, and the thin rings constructed by Fraenkel \cite{Fra1,Fra2}. On the basis of his numerical computations, Norbury \cite{Nor73} proposed a single family connecting Hill's vortex to the thin rings. Figure~\ref{fig1-1} illustrates the meridional sections of these rings.

Amick and Turner \cite{AT88} constructed a connected, unbounded set of solutions starting from Hill's vortex in the product of the flux axis and an energy space. It is natural to ask whether this continuum reaches arbitrarily large flux and whether the rings become thin as the flux tends to infinity. To answer the first question, we estimate the energy on bounded flux intervals. We also need this bound to apply the uniqueness theorem of Amick and Fraenkel for small flux \cite{AF88}. Cao et al.~\cite{CQZZ} showed that, as the impulse approaches that of Hill's vortex, the centered variational maximizers converge to Hill's vortex and their flux tends to zero. The local uniqueness theorem then identifies these rings with Norbury's local family.

In this paper, we bound the energy and the support of every centered finite-energy solution with bounded flux. It follows that the Amick--Turner continuum reaches every nonnegative flux and that the local uniqueness theorem near Hill's vortex applies to every centered finite-energy solution at sufficiently small flux. We then describe the shape of the entire core as the flux tends to infinity and prove uniqueness at sufficiently large flux.

\subsection{Main results}

We consider solutions of the form
\begin{equation*}
	u(t,x)=U(x-Wt e_z),\qquad U(x)\to0\quad\text{as }|x|\to\infty,
\end{equation*}
where $W>0$ is the traveling speed. In cylindrical coordinates $(r,\theta,z)$, we write
\begin{equation*}
	U=U^r(r,z)e_r+U^z(r,z)e_z,\qquad
	\nabla\times U=\omega^\theta(r,z)e_\theta .
\end{equation*}
By incompressibility there is a stream function $\psi$, vanishing on the axis, such that
\begin{equation*}
	U^r=-\frac1r\partial_z\psi,\qquad U^z=\frac1r\partial_r\psi.
\end{equation*}
The stream function for the relative velocity $U-We_z$ is $\psi-Wr^2/2$. By the steady vorticity equation, $\omega^\theta/r$ is constant along the streamlines of $U-We_z$. We take
\begin{equation*}
	\frac{\omega^\theta}{r}=\lambda\boldsymbol 1_D,\qquad
	D:=\left\{\psi-\frac W2r^2>k\right\},\qquad \lambda>0,\quad k\ge0.
\end{equation*}
Here $D$ is the meridional section of the vortex core for $r>0$, and we call the boundary value $k$ of the relative stream function the flux. We write $|E|$ for the Lebesgue measure of a set $E$ in its ambient Euclidean space. In estimates, $c$ and $C$ denote positive constants that may change from line to line. We indicate their dependence on parameters when needed.

By scaling we may take $W=2$ and $\lambda=15$. We set $\psi=r^2v$ and regard $v$ as a function of $x=(x',z)\in\R^4\times\R$, with $r=|x'|$. For $r>0$, the inequality defining $D$ becomes $F>0$, where
\begin{equation*}
	F(r,z):=v(r,z)-1-\frac{k}{r^2}.
\end{equation*}
For $k=0$, we set $F=v-1$ also on the axis and take $D=\{F>0\}$ in the closed meridional half-plane. For $k>0$, we extend $F_+$ by zero to the axis. The corresponding set in $\R^5$ is
\begin{equation*}
	D^{(5)}:=\{(x',z)\in\R^4\times\R:(|x'|,z)\in D\}.
\end{equation*}
Since $D$ is the positive set of $F$, we have $\|F_+\|_{L^\infty(\R^5)}=\|F\|_{L^\infty(D)}$. The identity
\begin{equation*}
	\frac{\omega^\theta}{r}
	=-\left(v_{rr}+\frac3r v_r+v_{zz}\right)
	=-\Delta_{\R^5}v,
\end{equation*}
leads to the free boundary problem
\begin{equation}\label{1-1}
	\left\{
	\begin{aligned}
		-\Delta_{\R^5}v&=15\boldsymbol 1_{D^{(5)}}\quad\text{in }\R^5,\\
		v(x)&\to0\quad\text{as }|x|\to\infty.
	\end{aligned}
	\right.
\end{equation}
For this stream function formulation see \cite{BF1,Ni,AF,AT88}. For the weak formulation in $\R^5$ see \cite{AF88}. In \eqref{1-1} we fix the speed and the value of $\omega^\theta/r$ in the core and vary the flux $k$.

A finite-energy solution is a nonzero weak solution $v=v(r,z)$ of \eqref{1-1} which tends to zero uniformly at infinity and satisfies
\begin{equation}\label{1-2}
	v\in D^{1,2}(\R^5)\cap W^{2,p}_{\mathrm{loc}}(\R^5)
	\quad\text{for every }1<p<\infty,\qquad
	\int_{\R^5}|\nabla v|^2\,dx<\infty.
\end{equation}
Here $D^{1,2}(\R^5)$ is the completion of $C_c^\infty(\R^5)$ in the gradient norm. The Dirichlet integral of $v$ is proportional to the kinetic energy of $U$. The local regularity in \eqref{1-2} follows from the boundedness of the right-hand side in \eqref{1-1}. We consider centered solutions in the sense of Amick and Fraenkel \cite{AF88}, for which
\begin{equation*}
	v(r,z)=v(r,-z),\qquad v_z(r,z)<0\quad(z>0).
\end{equation*}
The section at fixed $r$ is an interval symmetric about $z=0$ whenever it is nonempty.

In this normalization, Hill's vortex occurs at $k=0$ and has the form
\begin{equation*}
	v_H(r,z)=
	\begin{cases}
		\frac{5-3(r^2+z^2)}{2}, & r^2+z^2\le 1,\\[4pt]
		(r^2+z^2)^{-3/2}, & r^2+z^2\ge 1.
	\end{cases}
\end{equation*}
Since the radial derivatives of the two formulas agree on $r^2+z^2=1$, we have $-\Delta_{\R^5}v_H=15\boldsymbol 1_{\{r^2+z^2<1\}}$ in $\R^5$. The cross-section is $D=\{(r,z): r\ge0,\ r^2+z^2<1\}$. Inside this core the relative stream function is $\psi-r^2=\frac{3}{2}r^2(1-r^2-z^2)$, whose small positive level curves lie off the axis. The nearby Norbury rings have an inner radius of order $\sqrt{k}$, and their cross-sections approach the half-disk as $k$ tends to zero \cite{Nor72} (see also \cite{AF88}).

We first estimate the energy and the size of the support uniformly when the flux is bounded.

\begin{theorem}\label{thm1-1}
	For every $K<\infty$ and $0<\alpha<1$ there is $C_{K,\alpha}$ such that
	\begin{equation*}
		\|v\|_{C^{1,\alpha}(\R^5)} +\|\nabla v\|_{L^2(\R^5)}\le C_{K,\alpha}
	\end{equation*}
	for every centered finite-energy solution of \eqref{1-1} with $0\le k\le K$. Moreover, the support satisfies
	\begin{equation*}
		D^{(5)}\subset B_{C_K}(0)\subset\R^5 .
	\end{equation*}
\end{theorem}

By Theorem~\ref{thm1-1}, the solutions have a common energy bound for $k$ in a bounded interval. We can therefore apply the uniqueness theorem of Amick and Fraenkel for small positive $k$. Using also the uniqueness of Hill's vortex \cite{AF}, we find exactly one centered finite-energy solution for each sufficiently small $k\ge0$. On the Amick--Turner continuum, boundedness of the flux would now imply boundedness in the energy space. Since the continuum is connected, its fluxes form an unbounded interval containing zero, so the interval is $[0,\infty)$. We refer to Corollaries~\ref{cor2-7} and \ref{cor2-8} for these two consequences.

For large flux we first show that the points in the core have radii close to $\sqrt{k/3}$. We then translate the center of the cross-section to the origin and scale by the radius of a disk of the same area.

\begin{samepage}
	\begin{theorem}[Large flux]\label{thm1-2}
		Let $(k,v)$ be a centered finite-energy solution of \eqref{1-1} with cross-section $D$. Define $R$, $R_{\mathrm{center}}$ and $a$ by
		\begin{equation*}
			R:=\sqrt{k/3},\quad
			R_{\mathrm{center}}:=\frac{1}{|D|}\int_D r\,dr\,dz,\qquad
			a:=\left(\frac{|D|}{\pi}\right)^{1/2}.
		\end{equation*}
		By symmetry, the average of $(r,z)$ over $D$ is $(R_{\mathrm{center}},0)$. We call this point the center of $D$. As $k\to\infty$,
		\begin{equation*}
			\frac{R_{\mathrm{center}}}{R}\to1,\qquad a\to0,\qquad
			a^2\ln\frac{R_{\mathrm{center}}}{a}\to\frac8{15}.
		\end{equation*}

		For all sufficiently large $k$, $D\subset B_{Ca}((R_{\mathrm{center}},0))$. As $k\to\infty$, the closure of $(D-(R_{\mathrm{center}},0))/a$ converges to the closed unit disk in Hausdorff distance. The boundary of the rescaled cross-section is a graph in polar coordinates, and its radius converges to one in $C^{1,\alpha}(\mathbb T)$ for every $0<\alpha<1$, where $\mathbb T=\R/(2\pi\mathbb Z)$. In particular, $D$ is simply connected and its revolution is a solid torus for all sufficiently large $k$. The bounds and limits are uniform over centered solutions.

		Moreover, there is $0<K_\infty<\infty$ such that \eqref{1-1} has exactly one centered finite-energy solution for each $k\ge K_\infty$. Writing $v_k$ for this solution, the map $k\mapsto v_k$ is continuous from $[K_\infty,\infty)$ into $D^{1,2}(\R^5)\cap C^{1,\beta}(\R^5)$ for every $0<\beta<1$.
	\end{theorem}
\end{samepage}

\begin{figure}[htbp]
	\centering
	\includegraphics[width=0.92\linewidth]{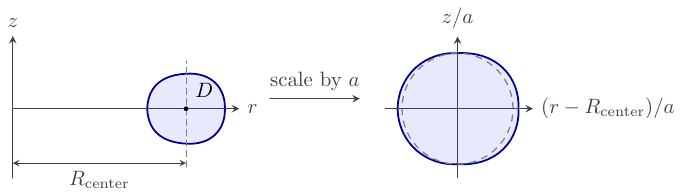}
	\caption{The cross-section at large flux, centered at $(R_{\mathrm{center}},0)$ and scaled by $a$. The solid curve on the right shows the boundary after scaling. The dashed circle is $\partial B_1(0)$.}
	\label{fig1-2}
\end{figure}

The centered solutions on the Amick--Turner continuum therefore agree with Norbury's local family near Hill's vortex and with the unique solution in Theorem~\ref{thm1-2} for each sufficiently large $k$. For $k_0<k<K_\infty$, with $k_0$ and $K_\infty$ as in Corollary~\ref{cor2-7} and Theorem~\ref{thm1-2}, it remains open whether the solution is unique and whether every cross-section is simply connected.

\begin{figure}[htbp]
	\centering
	\includegraphics[width=0.92\linewidth]{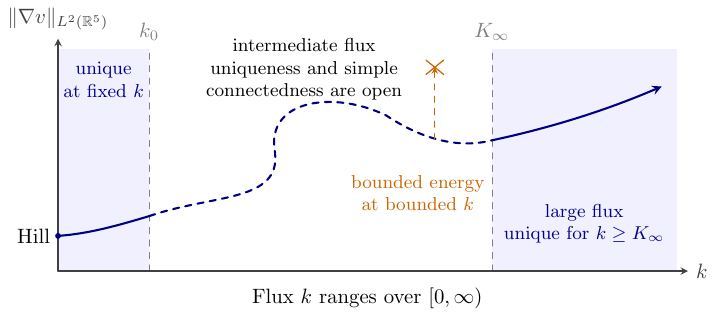}
	\caption{The Amick--Turner continuum covers every nonnegative flux. The solid parts correspond to the intervals of flux on which we prove uniqueness. The dashed part is a sketch of the continuum between these intervals.}
	\label{fig1-3}
\end{figure}

\subsection{Related works}

Elliptic constructions of vortex rings go back to Fraenkel and Berger \cite{BF1} and Ni \cite{Ni}. Friedman and Turkington \cite{FT} used Benjamin's variational formulation \cite{Ben} to obtain existence and concentration results. For other variational constructions, including mountain-pass and rearrangement methods, see \cite{AS,Bu87a,Bu,Bu03,Bad}. Buffoni \cite{Buf} constructed nested axisymmetric rings, and Abe \cite{Abe22} constructed rings with swirl in Beltrami flows.

Following Fraenkel's construction of vortex rings of small cross-section, there have been further results on desingularization and asymptotics of concentrated rings \cite{BF80,YJ,DV}.

Cao et al.~\cite{CWZ} desingularized a circular filament to construct vortex patches of small cross-section. Cao et al.~\cite{CWWZ} constructed concentrated vortex rings in the whole space for more general vorticity functions and proved concentration of the cross-sections. Cao and Zhan \cite{CZ} constructed thin vortex rings with and without swirl in the whole space. See Cao et al.~\cite{CWZ19} for vortex rings with swirl in an infinite pipe and outside a ball. Ao et al.~\cite{Ao} constructed several nearby vortex rings traveling together.

Friedman and Turkington \cite{FT} proved that the rescaled cross-sections of their constrained maximizers converge to a disk. Cao et al.~\cite{CQYZZ} obtained the Rankine potential as the limit of the rescaled potentials for vortex rings of small cross-section. Cao et al.~\cite{CLQZZ} obtained a radial limit when the vorticity function is a power on the positive half-line and vanishes elsewhere. We show that for every centered solution the entire core lies close to a circle and the rescaled cross-section tends to a disk as the flux tends to infinity.

Since the circulation is $15\pi R_{\mathrm{center}}a^2$, the leading term of the Kelvin--Hicks speed formula \cite{Fra1} (see also Cao et al.~\cite{CQYZZ}) is $\frac{15}{4}a^2\ln(R_{\mathrm{center}}/a)$. By Theorem~\ref{thm1-2}, this expression tends to the prescribed speed $W=2$. For a variational derivation of the Kelvin--Hicks formula, see Abe and Woo \cite{AW26}.

Cao et al.~\cite{CQYZZ} proved uniqueness for thin toroidal patches away from the axis with the same vortex strength, circulation and traveling speed. Their comparison uses local Pohozaev identities. Cao et al.~\cite{CLQZZ} proved uniqueness for a family of thin vortex rings whose vorticity function is a power on the positive half-line and vanishes elsewhere. We prescribe the flux $k$, so the circulations $15\int_D r\,dr\,dz$ of the two solutions may differ.

The proof in Section~\ref{sec5} follows the author's contour dynamics argument for concentrated rotating patches and helical vortices \cite{Qin}. We write the two rescaled boundaries as graphs in polar coordinates and subtract their boundary equations. The linearization at a disk controls the higher Fourier modes of the difference. To estimate the difference of the cross-sectional areas, we retain the constant terms in the equations for the potentials. We use the radial integral identity to bound the difference of the two values of $R_{\mathrm{center}}$.

For the response of Hill's vortex to axisymmetric perturbations, see \cite{MM,Poz,PE}. Choi \cite{Choi20} proved its orbital stability under axisymmetric perturbations, and Choi and Jeong \cite{CJ} proved that filamentation occurs for perturbations of Hill's vortex. Cao et al.~\cite{CQZZ} proved orbital stability of a Norbury ring close to Hill's vortex under axisymmetric perturbations. For nonlinear stability of vortex rings of small cross-section, see \cite{CQYZZ,CLQZZ}. Protas \cite{Pro} computed the spectrum numerically for axisymmetric perturbations of Norbury rings and found neutral stability for thin rings and instability for rings of sufficiently large cross-section. Understanding the stability of solutions along the global continuum remains an open problem.

For the vortex atmosphere, namely the fluid that lies outside the core of a traveling vortex and moves with it, see Choi et al.~\cite{CJS2}. The evolution of concentrated vortex rings has also been studied \cite{BCM,But2,CJ25,DDMW,DHLM,EY,Fen,GS0,GS2,GHH,GJZ,GMT}.

The variational family of traveling dipoles obtained by Abe et al.~\cite{AJQZ} is a two-dimensional analogue of the Norbury family. Uniqueness is known near the Lamb dipole and for large impulse, and remains open for intermediate values. Abe and Choi \cite{Abe} proved orbital stability of the Lamb dipole. Under symmetry and sign assumptions, Abe et al.~\cite{ACJ} proved stability for perturbations that may have infinite mass. See \cite{AJY,CJY} for stability of several Lamb dipoles and of vortex quadrupoles with odd symmetry in both variables, respectively.

A related question about the end of a global branch concerns translating patches that meet at a point, as in Sadovskii's vortex \cite{Sad}. Early studies of translating patch pairs and touching configurations include \cite{Nor75,Pie,ST}. For the variational construction of planar vortices we refer to \cite{T83a,T83b}. Hmidi and Mateu \cite{HM} constructed local branches from pairs of point vortices, and Garc\'ia and Haziot \cite{GH} continued these branches globally. Cao et al.~\cite{CQZZ2} proved uniqueness and orbital stability for concentrated symmetric traveling patch pairs. For existence, variational stability and nondegeneracy results concerning touching vortex pairs, see \cite{CJS,HT,ACJSW,DMW}. It remains open whether the branch of counter-rotating pairs in \cite{GH} converges to a pair of patches that touch.

A similar question arises for two-dimensional rotating vortex patches. Hassainia et al.~\cite{HasMW} proved global continuation of the branches bifurcating from the disk, extending the local results in \cite{Burb,Hmi}. For other results on the construction and symmetry of rotating patches, see \cite{Cas4,GPSY,CLGS}. The numerical computations of Wu et al.~\cite{WOZ} suggested the formation of corners in the limit (see \cite{Ove,WZZ26} for the local analysis). Huang et al.~\cite{HTZ} constructed $m$-fold symmetric rotating patches with right-angle corners for every $m\ge12$. The conjecture in \cite{HasMW} that patches with such corners arise as limits along the branches from the disk remains open. For Norbury rings, we show that the continuum from Hill's vortex reaches every nonnegative flux and that the rescaled cross-sections tend to a disk as the flux tends to infinity.

\subsection{Outline of the proof}

In Section~\ref{sec2}, we first estimate $F_+$ and the gradient of $v$. To control the height of the cross-section, we derive the boundary inequality \eqref{2-19} from the equation for $-v_z$ and its positivity in $\{z>0\}$. If the height is large compared with $\sqrt{\|F_+\|_{L^\infty(\R^5)}}$, we rescale at a point where the axial derivative is small. The limit is independent of the axial variable. We choose a test function near the right endpoint of the interval on which the limit is positive. Its energy is strictly smaller than the lower bound required by \eqref{2-19}. Combining the resulting height bound with the Pohozaev identity, we obtain uniform volume and energy estimates. We then bound $r$ on $D$ by Sobolev's inequality and the gradient estimate.

In Section~\ref{sec3}, we show that $D$ lies near $r=R$ when the flux is large. For the kernel estimates we use the expansions of Gallay and \v{S}ver\'ak \cite{GS15}. Integrating $\partial_r(r^{3/2}F)$ in $D$ on either side of $r=T$, we obtain
\begin{align*}
	T^{3/2}\int_\R F_+(T,z)\,dz
	&=\int_{D\cap\{r<T\}}\partial_r(r^{3/2}F)\,dr\,dz\\
	&=-\int_{D\cap\{r>T\}}\partial_r(r^{3/2}F)\,dr\,dz.
\end{align*}
Write $H(T)$ for half the length of the section at $r=T$ and $S(T)$ for the area of $D$ on the side under consideration. Since $F$ vanishes at the endpoints of the section, the left-hand side is at most $T^{3/2}\|F_z\|_{L^\infty(D)}H(T)^2$. If the smallest or largest radius of $D$ differs from $R$ by a fixed fraction of $R$, we bound the corresponding integral over $D\cap\{r<T\}$ or $D\cap\{r>T\}$ from below by $c\sqrt R\,S(T)$ for $T$ in an interval of length comparable to $R$. By $\|F_z\|_{L^\infty(D)}\le C\sqrt R$, this implies $H(T)^2\ge cR^{-3/2}S(T)$. Since $|S'(T)|=2H(T)$, integrating $|(\sqrt S)'|=H/\sqrt S$ shows that the portion of $D$ between the endpoints of the interval has area at least $cR^{1/2}$. Its mass is at least $cR^2$, contradicting the upper bound $CR^{3/2}$.

In Section~\ref{sec4} we use the same identity to bound the radial width of $D$ by $Ca$. Let $\ell$ be the larger of the distances from $R_{\mathrm{center}}$ to the smallest and largest radii of $D$. By estimating the integral over pairs of points on opposite sides of $r=T$, we have $H(T)^2\ge c(a/\ell)S(T)$ for $T$ in an interval of length $\ell/4$. The area of $D$ between the two lines is therefore at least $ca\ell$ and at most $\pi a^2$. Once $D$ lies within $Ca$ of its center, we evaluate the boundary equation to determine the limit of $a^2\ln R$.

We derive the equation for the limiting cross-section as in \cite{Qin}. The limiting potential is the logarithmic potential of the limiting cross-section plus a constant. It vanishes on the boundary of that cross-section. Thus the logarithmic potential is constant there, and Fraenkel's characterization \cite{Fb} (see also \cite{Rei}) identifies the limiting cross-section as a disk. As in \cite{Qin}, the rescaled potential has opposite signs on the two circles of an annulus about the unit circle and decreases strictly along each ray between them. The boundary lies in this annulus and meets each ray at the unique zero there.

Finally, in Section~\ref{sec5} we compare two cross-sections at the same flux. We use the radial integral identity to bound the difference of the centers. On the $m$th cosine mode, the linearized equation multiplies the difference of the radii by $1-1/m$. This controls the modes of order $m\ge2$. The difference of the areas and the vanishing first moments about the centers give the estimates for the mean and the first cosine mode, respectively. Combining these estimates, we obtain uniqueness.

\section{Compactness at bounded flux}\label{sec2}

Amick and Fraenkel \cite{AF88} bounded the support when both the flux and the energy are bounded. We first bound the supremum of $F_+$ and the height of the core, and then use a Pohozaev identity to deduce the energy and support bounds in Theorem~\ref{thm1-1}. These estimates yield compactness for bounded flux. They also show that $k$ is unbounded on the Amick--Turner continuum.

\subsection{Bounds for the potential and its gradient}\label{sec2-1}

We first bound $v$ and its gradient in terms of $k$. We use these bounds for large flux in Section~\ref{sec3}.

\begin{lemma}\label{lem2-1}
	There is a universal $C$ such that every centered finite-energy solution $(k,v)$ of \eqref{1-1} satisfies
	\begin{equation}\label{2-1}
		\|v\|_{L^\infty(\R^5)}\le C(1+\sqrt k),\qquad
		\|\nabla v\|_{L^\infty(\R^5)}\le C(1+\sqrt k)^{1/2}.
	\end{equation}
	If $0\le k\le K$ and $0<\alpha<1$, then
	\begin{equation*}
		\|v\|_{C^{1,\alpha}(\R^5)}\le C_{K,\alpha}.
	\end{equation*}
\end{lemma}

\begin{proof}
	The decay of $v$ implies that $D^{(5)}\subset\{v>1\}$ is bounded. The difference of $v$ and the Newtonian potential of $15\boldsymbol 1_{D^{(5)}}$ is harmonic on $\R^5$ and vanishes at infinity. Thus
	\begin{equation*}
		v(x)=\frac{15}{8\pi^2}\int_{D^{(5)}}|x-y|^{-3}\,dy\ge0.
	\end{equation*}
	Let $M=v(x_0)=\|v\|_{L^\infty(\R^5)}>0$. Using $|\mathbb S^4|=8\pi^2/3$,
	\begin{equation*}
		\frac{15}{8\pi^2}\int_{B_{\sqrt M/2}(x_0)}|x_0-y|^{-3}\,dy
		=5\int_0^{\sqrt M/2}t\,dt=\frac58M.
	\end{equation*}
	If $x\in B_{\sqrt M}(x_0)$ and $|y-x_0|\ge\sqrt M/2$, then $|x-y|\le3|x_0-y|$. Thus
	\begin{equation*}
		v(x)\ge\frac1{27}\frac{15}{8\pi^2}
		\int_{D^{(5)}\setminus B_{\sqrt M/2}(x_0)}|x_0-y|^{-3}\,dy
		\ge\frac{M-(5/8)M}{27}\ge cM .
	\end{equation*}
	Suppose $M>2/c$, where $c$ is the constant in the last inequality. For $x\in B_{\sqrt M}(x_0)$ with $r^2>2k/(cM)$, we have
	\begin{equation*}
		v-1-\frac{k}{r^2}\ge cM-1-\frac{cM}{2}
		=\frac{cM}{2}-1>0 .
	\end{equation*}
	Thus $B_{\sqrt M}(x_0)\cap\{r^2>2k/(cM)\}\subset D^{(5)}$.

	To bound the integral of the Newtonian kernel over $\{r^2\le2k/(cM)\}$, write $x=(x',z)$ and $x_0=(x_0',z_0)$. For $x'\neq x_0'$,
	\begin{equation*}
		\int_{\R}\frac{dz}{(|x'-x_0'|^2+(z-z_0)^2)^{3/2}}
		=\left[\frac{z-z_0}{|x'-x_0'|^2\sqrt{|x'-x_0'|^2+(z-z_0)^2}}\right]_{-\infty}^{\infty}
		=\frac{2}{|x'-x_0'|^2}.
	\end{equation*}
	Integrating first in $z$ and using that $|x'|^{-2}$ is radially decreasing in $\R^4$, we obtain
	\begin{align*}
		\int_{\{r^2<2k/(cM)\}}|x_0-x|^{-3}\,dx
		&=2\int_{|x'|^2<2k/(cM)}|x'-x_0'|^{-2}\,dx'\\
		&\le2\int_{|x'|^2<2k/(cM)}|x'|^{-2}\,dx'\\
		&=2|\mathbb S^3|\int_0^{\sqrt{2k/(cM)}} t\,dt
		=\frac{2|\mathbb S^3|k}{cM}\le \frac{Ck}{M}.
	\end{align*}
	For $M>2/c$, the integral over $D^{(5)}$ is at least the integral on $B_{\sqrt M}(x_0)$ minus the integral over $\{r^2\le2k/(cM)\}$. Thus
	\begin{align*}
		M
		&\ge\frac{15}{8\pi^2}\int_{B_{\sqrt M}(x_0)\cap\{r^2>2k/(cM)\}}|x_0-y|^{-3}\,dy\\
		&\ge 5\int_0^{\sqrt M}t\,dt-\frac{Ck}{M}
		=\frac52M-\frac{Ck}{M}.
	\end{align*}
	Hence $M\le2/c$ or $\frac32 M^2\le Ck$, which proves the bound for $v$ in \eqref{2-1}.

	By the interior gradient estimate on balls of radius $\sqrt M$,
	\begin{equation*}
		\|\nabla v\|_{L^\infty(\R^5)}
		\le C(M/\sqrt M+15\sqrt M)\le C(1+\sqrt k)^{1/2}.
	\end{equation*}
	For $p>5/(1-\alpha)$, the interior $W^{2,p}$ estimate and Morrey's inequality imply
	\begin{equation}\label{2-14}
		\|v\|_{C^{1,\alpha}(B_{1/2}(x))}
		\le C_p\bigl(\|v\|_{L^\infty(B_1(x))}+\|\Delta_{\R^5}v\|_{L^p(B_1(x))}\bigr).
	\end{equation}
	By \eqref{1-1} and \eqref{2-1}, we bound the right-hand side uniformly in $x$. This bounds the H\"older quotient for pairs with $|x-y|<1/2$. For $|x-y|\ge1/2$,
	\begin{equation*}
		\frac{|\nabla v(x)-\nabla v(y)|}{|x-y|^\alpha}
		\le2^{1+\alpha}\|\nabla v\|_{L^\infty(\R^5)} .
	\end{equation*}
	Together the two bounds control $\|v\|_{C^{1,\alpha}(\R^5)}$.
\end{proof}

\subsection{The vertical support estimate}\label{sec2-2}

In $\{z>0\}$, the $z$ derivative of the characteristic function in \eqref{1-1} is a measure supported on the free boundary. Since $-v_z>0$ for $z>0$, testing the differentiated equation with $\phi^2/(-v_z)$ leads to the following boundary estimate.

\begin{lemma}\label{lem2-2}
	Let $(k,v)$ be a centered finite-energy solution of \eqref{1-1}. For every compactly supported axisymmetric Lipschitz function $\phi$ in $\{z>0\}$,
	\begin{equation}\label{2-19}
		15\int_{\partial D\cap\{z>0\}}\frac{r^3\phi^2}{|\nabla F|}\,ds
		\le\int_0^\infty\int_0^\infty
		r^3(\phi_r^2+\phi_z^2)\,dr\,dz.
	\end{equation}
	Here $ds$ denotes the element of arc length in the $(r,z)$-plane.
\end{lemma}

\begin{proof}
	By the local regularity of $v$ and its strict decrease in $z>0$, the function $w=-v_z$ satisfies
	\begin{equation*}
		w\in W^{1,p}_{\mathrm{loc}}(\R^4\times(0,\infty))\cap C(\R^4\times(0,\infty)),\qquad w>0
		\qquad(p>5).
	\end{equation*}
	By evenness and $F_z=-w<0$ for $z>0$, a nonempty vertical section at $r>0$ is an interval $(-H(r),H(r))$. At its upper endpoint the implicit function theorem shows that the boundary is a $C^{1,\alpha}$ graph $z=H(r)$, and
	\begin{equation*}
		H'(r)=-\frac{F_r}{F_z},\qquad
		ds=\sqrt{1+H'(r)^2}\,dr=\frac{|\nabla F|}{w}\,dr.
	\end{equation*}
	Let $\eta$ be a smooth axisymmetric test function compactly supported in $z>0$. We differentiate the weak equation for $v$ in $z$, divide by $|\mathbb S^3|$, and use $\int_0^{H(r)}\eta_z(r,z)\,dz=\eta(r,H(r))$ for each nonempty section. The weight $r^3$ vanishes on the axis when the core meets it. Thus
	\begin{equation}\label{2-16}
		\begin{split}
			\int_{r,z>0}r^3\nabla w\cdot\nabla\eta\,dr\,dz &=15\int_Dr^3\eta_z\,dr\,dz\\
			&=15\int_{\partial D\cap\{z>0\}} \frac{r^3w\eta}{|\nabla F|}\,ds
		\end{split}
	\end{equation}
	Since $w>0$ on $\operatorname{supp}\phi$, we have $\phi^2/w\in W^{1,p}_c(\R^4\times(0,\infty))$. Using a radial mollifier, we approximate $\phi^2/w$ by smooth axisymmetric functions $\eta_j$ supported in a fixed compact subset of $\{z>0\}$. Then $\eta_j\to\phi^2/w$ in $W^{1,p}(\R^4\times(0,\infty))$ and uniformly on $\R^4\times(0,\infty)$. Choose $C_\phi$ so that the supports of $\eta_j$ and $\phi^2/w$ lie in $r<C_\phi$. On this region the graph formula implies
	\begin{equation*}
		\int_{\partial D\cap\{0<r<C_\phi,\,z>0\}}
		\frac{r^3w}{|\nabla F|}\,ds
		\le\int_0^{C_\phi}r^3\,dr=\frac{C_\phi^4}{4}.
	\end{equation*}
	H\"older's inequality for the interior integral and the preceding bound for the boundary integral give
	\begin{align*}
		\left|\int_{r,z>0} r^3\nabla w\cdot\nabla(\eta_j-\phi^2/w)\,dr\,dz\right|
		&\le C_{\phi,w}\|\eta_j-\phi^2/w\|_{W^{1,p}(\R^4\times(0,\infty))}\to0,\\
		15\left|\int_{\partial D\cap\{z>0\}}\frac{r^3w(\eta_j-\phi^2/w)}{|\nabla F|}\,ds\right|
		&\le C_{\phi}\|\eta_j-\phi^2/w\|_{L^\infty(\R^4\times(0,\infty))}\to0.
	\end{align*}
	We may thus test \eqref{2-16} with $\phi^2/w$. Completing the square, we obtain
	\begin{equation}\label{2-18}
		\begin{split}
			15\int_{\partial D\cap\{z>0\}} \frac{r^3\phi^2}{|\nabla F|}\,ds
			&=\int_{r,z>0}r^3\left( \frac{2\phi}{w}\nabla w\cdot\nabla\phi
			-\frac{\phi^2}{w^2}|\nabla w|^2\right)\,dr\,dz\\
			&=\int_{r,z>0}r^3|\nabla\phi|^2\,dr\,dz -\int_{r,z>0}r^3 \left|\nabla\phi-\frac\phi w\nabla w\right|^2\,dr\,dz\\
			&\le\int_{r,z>0}r^3|\nabla\phi|^2\,dr\,dz.
		\end{split}
	\end{equation}
\end{proof}

Before rescaling, we choose a point away from the axis where $F>0$ and where the derivative in $z$ will be small after rescaling. We also need a bound on $F_+$ for $z$ in an interval about this point.

\begin{lemma}\label{lem2-3}
	Let $(k,v)$ be a centered finite-energy solution of \eqref{1-1}, and suppose that for some $N\ge1$,
	\begin{equation*}
		\sup_{(r,z)\in D}|z|>2N\sqrt{\|F_+\|_{L^\infty(\R^5)}} .
	\end{equation*}
	There are $\mu>0$ and $(r_0,z_0)$ such that
	\begin{align*}
		r_0&\ge c\sqrt{\mu},\qquad z_0>\frac{N\sqrt{\mu}}4,\\
		F(r_0,z_0)&\ge\frac{\mu}8,\qquad 0<-F_z(r_0,z_0)\le\frac{C\sqrt{\mu}}N ,
	\end{align*}
	where $c,C>0$ are universal. Moreover,
	\begin{equation}\label{2-6}
		0\le F_+(r,z)\le4\mu
		\qquad\left(r\ge0,\ |z-z_0|<\frac{N\sqrt{\mu}}4\right).
	\end{equation}
\end{lemma}

\begin{proof}
	Let $M(z):=\sup_{r>0}F_+(r,z)$. By symmetry and monotonicity, $M(0)=\|F_+\|_{L^\infty(\R^5)}$ and $M$ is nonincreasing on $[0,\infty)$. By the height assumption and symmetry, there is a point of $D$ with $z>2N\sqrt{M(0)}$. If the point is on the axis, then $k=0$ and by continuity we can find a point with positive radius at the same height. Thus $M(2N\sqrt{M(0)})>0$ by monotonicity. By Lemma~\ref{lem2-1},
	\begin{equation*}
		|M(z)-M(z')|\le\|v_z\|_{L^\infty(\R^5)}|z-z'|.
	\end{equation*}
	Choose $t$ maximizing $z+N\sqrt{M(z)}$ over $[0,2N\sqrt{M(0)}]$, and set $\mu:=M(t)>0$. Comparing with the right endpoint, we have
	\begin{equation*}
		t+N\sqrt{\mu}\ge2N\sqrt{M(0)}+N\sqrt{M(2N\sqrt{M(0)})}>2N\sqrt{M(0)} .
	\end{equation*}
	Since $\mu\le M(0)$, it follows that $t>N\sqrt{M(0)}\ge N\sqrt{\mu}$. Comparing with $z=t-N\sqrt{\mu}$, we have
	\begin{equation*}
		t-N\sqrt{\mu}+N\sqrt{M(t-N\sqrt{\mu})}\le t+N\sqrt{\mu}.
	\end{equation*}
	Hence $M(t-N\sqrt{\mu})\le4\mu$, and $F_+(r,z)\le4\mu$ whenever $z\ge t-N\sqrt{\mu}$.

	To find a point away from the axis, choose $r>0$ with $F(r,t)\ge\mu/2$ and set $x_0=(r e_1,t)$, where $e_1$ is the first unit vector in $\R^4$. On $\{r>0\}$, $r^{-2}$ is harmonic. By \eqref{1-1} and the distributional inequality $-\Delta_{\R^5}F_+\le-\boldsymbol 1_{\{F>0\}}\Delta_{\R^5}F$, we have $-\Delta_{\R^5}F_+\le15$ there. For $k>0$, $F_+$ vanishes near the axis, and for $k=0$ the inequality holds on $\R^5$ by applying it to $v-1$. Since $\Delta_{\R^5}|x-x_0|^2=10$, the function $F_+(x)+\frac32|x-x_0|^2$ is subharmonic. Its mean on $B_{\sqrt\mu/4}(x_0)$ is at least its value at $x_0$, so
	\begin{equation*}
		\frac1{|B_{\sqrt{\mu}/4}(x_0)|}\int_{B_{\sqrt{\mu}/4}(x_0)}F_+\,dx
		\ge\frac{\mu}{2}-\frac32\frac{\mu}{16}
		=\frac{13}{32}\mu>\frac{\mu}{4}.
	\end{equation*}
	The intersection of this ball with $\{|x'|<c\sqrt{\mu}\}$ has length at most $\sqrt{\mu}/2$ in the $z$ direction. Its sections at fixed $z$ have volume at most $(\pi^2/2)c^4\mu^2$. Dividing the volume of the intersection by that of the ball, we obtain the bound
	\begin{equation*}
		\frac{(\sqrt{\mu}/2)(\pi^2/2)c^4\mu^2}{(8\pi^2/15)(\sqrt{\mu}/4)^5}\le Cc^4 .
	\end{equation*}
	The ball lies above $t-N\sqrt{\mu}$ since $N\ge1$. If $F_+<\mu/8$ everywhere in the ball with $|x'|\ge c\sqrt{\mu}$, its mean would be at most $\mu/8+4Cc^4\mu<\mu/4$ for $c$ small enough. Thus there is a point of radius $r_0\ge c\sqrt{\mu}$ and height at least $t-\sqrt{\mu}/4$ with $F\ge\mu/8$.

	The interval $[t-3N\sqrt{\mu}/4,t-N\sqrt{\mu}/2]$ lies between $t-N\sqrt{\mu}$ and $t-\sqrt{\mu}/4$. The upper bound for $F_+$ and monotonicity of $F$ in $z>0$ therefore imply
	\begin{equation*}
		\frac{\mu}{8}\le F(r_0,z)\le4\mu
		\qquad \left(t-\frac{3N\sqrt{\mu}}4\le z\le t-\frac{N\sqrt{\mu}}2\right).
	\end{equation*}
	By the mean value theorem, there is $z_0$ in this interval such that
	\begin{align*}
		0<-F_z(r_0,z_0)
		&=\frac{F(r_0,t-3N\sqrt{\mu}/4)-F(r_0,t-N\sqrt{\mu}/2)}{N\sqrt{\mu}/4}
		\le\frac{C\sqrt{\mu}}N .
	\end{align*}
	Finally, $z_0-N\sqrt{\mu}/4\ge t-N\sqrt{\mu}>0$. The bound for $F_+$ above $t-N\sqrt{\mu}$ proves \eqref{2-6}.
\end{proof}

In the rescaled variables the distance from the origin to the axis may have a finite or an infinite limit. For the corresponding weights, we construct a test function whose energy is strictly smaller than the boundary term at the right endpoint of the interval on which the limit is positive.

\begin{lemma}\label{lem2-4}
	Let $0<r_*\le\infty$, and set
	\begin{equation*}
		p(t)=
		\begin{cases}
			(1+t/r_*)_+^3,& r_*<\infty,\\
			1,&r_*=\infty.
		\end{cases}
		\qquad I:=\{t:p(t)>0\}.
	\end{equation*}
	Let $f\in C^1(I)$ be bounded above with $f(0)>0$. Suppose that $-(pf')'=15p$ on the component $(t_-,t_+)$ of $\{t\in I:f(t)>0\}$ containing $0$. Then $t_+\in I$, $f(t_+)=0$ and $f'(t_+)<0$. There is a compactly supported Lipschitz function $\varphi$ on $\R$ with $\varphi(t_+)=1$ and
	\begin{equation}\label{2-7}
		\int_I p(t)|\varphi'(t)|^2\,dt<\frac{15p(t_+)}{|f'(t_+)|}.
	\end{equation}
\end{lemma}

\begin{proof}
	Integrating $-(pf')'=15p$ from $0$ to $t\in(t_-,t_+)$ and using $p(0)=1$, we find
	\begin{align*}
		p(t)f'(t)&=f'(0)-15\int_0^t p(s)\,ds,\\
		f'(t)&=
		\begin{cases}
			-\dfrac{15}{4}(r_*+t)+\dfrac{r_*^3(f'(0)+15r_*/4)}{(r_*+t)^3},&r_*<\infty,\\[4pt]
			f'(0)-15t,&r_*=\infty .
		\end{cases}
	\end{align*}
	If $t_+=\infty$, in either case $f'(t)\to-\infty$, which contradicts $f>0$ on $(0,\infty)$. Thus $t_+$ is a finite point of $I$ and $f(t_+)=0$. By the mean value theorem, $f'<0$ at some point of $(0,t_+)$. Since $pf'$ is strictly decreasing there and $p(t_+)>0$, we have $f'(t_+)<0$.

	First suppose $r_*<\infty$. If $t_->-r_*$, then $f(t_-)=0$ and $f'(t_-)\ge0$. If $t_-=-r_*$, the term proportional to $(r_*+t)^{-3}$ in $f'$ would produce a term proportional to $(r_*+t)^{-2}$ in $f$. Since $f$ is positive and bounded above on $(t_-,t_+)$, the coefficient must vanish. In this case $\lim_{t\downarrow t_-}(r_*+t)^3f'(t)=0$. Setting $r_+:=r_*+t_+$, we obtain
	\begin{align*}
		r_+^3|f'(t_+)|
		&=15\int_{t_-}^{t_+}(r_*+t)^3\,dt
		-\lim_{t\downarrow t_-}(r_*+t)^3f'(t)\\
		&\le 15\int_0^{r_+}r^3\,dr=\frac{15}{4}r_+^4.
	\end{align*}
	Take $\varphi=1$ on $[-r_*,t_+]$ and let it decrease linearly to zero on $[t_+,t_++r_+]$, with a compactly supported Lipschitz extension to the left of $-r_*$. Using $r=r_*+t$, we have
	\begin{equation*}
		r_*^3\int_I p(t)|\varphi'(t)|^2\,dt
		=\frac1{r_+^2}\int_{r_+}^{2r_+}r^3\,dr=\frac{15}{4}r_+^2
		<4r_+^2\le \frac{15r_+^3}{|f'(t_+)|}.
	\end{equation*}
	This proves \eqref{2-7} because $r_*^3 p(t_+)=r_+^3$.

	If $r_*=\infty$, take $\varphi(t)=(1-|t-t_+|/\ell)_+$. Its energy is
	\begin{equation*}
		\int_{t_+-\ell}^{t_++\ell}\ell^{-2}\,dt=\frac2\ell .
	\end{equation*}
	For $\ell$ large this is strictly less than $15/|f'(t_+)|$.
\end{proof}

We now apply the preceding estimate to the limit obtained by rescaling at the points in Lemma~\ref{lem2-3}.

\begin{lemma}\label{lem2-5}
	There is a universal constant $C$ such that every centered finite-energy solution $(k,v)$ of \eqref{1-1} satisfies
	\begin{equation*}
		\sup_{(r,z)\in D}|z|\le C\sqrt{\|F_+\|_{L^\infty(\R^5)}}.
	\end{equation*}
\end{lemma}

\begin{proof}
	If the estimate fails, there are solutions $(k_n,v_n)$ with cross-sections $D_n$ such that
	\begin{equation*}
		\sup_{(r,z)\in D_n}|z|>2n\sqrt{\|(F_n)_+\|_{L^\infty(\R^5)}}.
	\end{equation*}
	By Lemma~\ref{lem2-3} with $N=n$, we can choose $(r_{0,n},z_n)$ and $\mu_n$. Rescaling about $(r_{0,n},z_n)$ by $\sqrt{\mu_n}$, we set
	\begin{align*}
		X&:=\left(\frac{r-r_{0,n}}{\sqrt{\mu_n}},\frac{z-z_n}{\sqrt{\mu_n}}\right),
		\qquad r_n:=\frac{r_{0,n}}{\sqrt{\mu_n}},\\
		\tilde F_n(X)&:=\frac{F_n(r_{0,n}+\sqrt{\mu_n}\,X_1,z_n+\sqrt{\mu_n}\,X_2)}{\mu_n}.
	\end{align*}
	Then $r_n\ge c$, and by \eqref{2-6}, $\tilde F_n\le4$ in $X_1>-r_n$, $|X_2|<n/4$. At the origin,
	\begin{align*}
		\tilde F_n(0)&=\frac{F_n(r_{0,n},z_n)}{\mu_n}\ge\frac18,\\
		0<-\partial_{X_2}\tilde F_n(0)&=\frac{-\partial_z F_n(r_{0,n},z_n)}{\sqrt{\mu_n}}\le\frac Cn.
	\end{align*}
	Set $p_n(t):=(1+t/r_n)_+^3$. The inequality $z_n>n\sqrt{\mu_n}/4$ ensures $z>0$ for $|X_2|<n/4$. For $X_1>-r_n$, we use that $r^{-2}$ is harmonic, write the equation for $F_n$ in the variables $X$, and differentiate in $X_2$. We obtain
	\begin{align*}
		-\nabla_X\cdot(p_n\nabla_X\tilde F_n)&=15p_n\boldsymbol 1_{\{\tilde F_n>0\}},\\
		-\nabla_X\cdot\bigl(p_n\nabla_X(-\partial_{X_2}\tilde F_n)\bigr)
		&=-15p_n\partial_{X_2}\boldsymbol 1_{\{\tilde F_n>0\}}\ge0.
	\end{align*}
	The last inequality holds in distributions since $\tilde F_n$ decreases in $X_2$ for $|X_2|<n/4$.

	Up to a subsequence, $r_n\to r_*\in[c,\infty]$. Then $p_n$ converges uniformly on compact intervals to the weight $p$ in Lemma~\ref{lem2-4}. On compact subsets of $I\times\R$, the weights also converge with their first derivatives and are bounded away from zero. After division by $p_n$, the equation for $4-\tilde F_n\ge0$ has coefficient $3/(r_n+X_1)$ in the first-order term and right-hand side bounded by $15$. The value of $4-\tilde F_n$ at the origin is at most $4$. Applying Harnack's inequality on a chain of overlapping balls starting from the origin, we bound $4-\tilde F_n$ on compact subsets of $I\times\R$. The interior estimates then bound $\|\tilde F_n\|_{W^{2,q}(B)}$ on each ball $B$ compactly contained in $I\times\R$, for every $1<q<\infty$. Passing to a further subsequence, we have
	\begin{align*}
		\tilde F_n&\to F_\infty\quad\text{in }C^1_{\mathrm{loc}}(I\times\R),\\
		F_\infty&\le4,\qquad F_\infty(0)\ge\frac18,\qquad
		-\partial_{X_2}F_\infty\ge0,\qquad \partial_{X_2}F_\infty(0)=0 .
	\end{align*}
	The interior bounds also give $F_\infty\in W^{2,q}_{\mathrm{loc}}(I\times\R)$. By the local convergence of $p_n$ and its derivative, $-\partial_{X_2}F_\infty$ is a weak supersolution of $-\nabla_X\cdot(p\nabla_X)$ on $I\times\R$. It is nonnegative and vanishes at the origin, so the strong maximum principle implies $F_\infty(X)=F_\infty(X_1)$. Let $(t_-,t_+)$ be the component of $\{t\in I:F_\infty(t)>0\}$ containing $0$. On compact subsets of $(t_-,t_+)\times\R$, $\tilde F_n>0$ for large $n$. Passing to the limit in the equation there, we have $-(pF_\infty')'=15p$ on $(t_-,t_+)$.

	By Lemma~\ref{lem2-4}, $F_\infty(t_+)=0$, $F_\infty'(t_+)<0$, and we can choose a test function $\varphi$ with $\varphi(t_+)=1$ satisfying \eqref{2-7}. Fix a nonzero $\eta\in C_c^\infty((-1,1))$ and use $\phi(X)=\varphi(X_1)\eta(X_2/L)$ for $L>0$. In the original variables this is an axisymmetric Lipschitz function with compact support in $z>0$ when $n>4L$. In \eqref{2-19}, the weight is $r^3=r_{0,n}^3p_n(X_1)$. The factors $\sqrt{\mu_n}$ from the arc length element and from $|\nabla F_n|$ cancel on the left. On the right, the change of area contributes $\mu_n$ and the squared gradient of the test function contributes $\mu_n^{-1}$. After cancellation we obtain
	\begin{equation}\label{2-8}
		15\int_{\partial\{\tilde F_n>0\}\cap\{X_1>-r_n,\,|X_2|<n/4\}}
		\frac{p_n(X_1)\phi^2}{|\nabla_X\tilde F_n|}\,ds
		\le \int_{-n/4}^{n/4}\int_{-r_n}^{\infty}p_n(X_1)|\nabla_X\phi|^2\,dX_1\,dX_2.
	\end{equation}
	By $C^1$ convergence, $\partial_{X_1}\tilde F_n<0$ near $X_1=t_+$ for $|X_2|\le L$, and $\tilde F_n$ has opposite signs at the two ends of a small interval about $t_+$. By the implicit function theorem, the zero set there is a graph $X_1=t_{+,n}(X_2)$ on $[-L,L]$ in the boundary of the positive set, and $t_{+,n}$ converges to $t_+$ in $C^1([-L,L])$. Since $t_{+,n}'=-\partial_{X_2}\tilde F_n/\partial_{X_1}\tilde F_n$ on the graph,
	\begin{equation*}
		\frac{ds}{|\nabla_X\tilde F_n|}=
		\frac{dX_2}{|\partial_{X_1}\tilde F_n(t_{+,n}(X_2),X_2)|}.
	\end{equation*}
	Restricting the boundary integral in \eqref{2-8} to this graph and integrating the two terms of $|\nabla_X\phi|^2$ separately in $X_1$ and $X_2$, we have
	\begin{align*}
		&15\int_{-L}^L
		\frac{p_n(t_{+,n}(X_2))\varphi(t_{+,n}(X_2))^2\,\eta(X_2/L)^2}
		{|\partial_{X_1}\tilde F_n(t_{+,n}(X_2),X_2)|}\,dX_2\\
		&\quad\le \int_\R p_n(t)|\varphi'(t)|^2\,dt\int_\R\eta(X_2/L)^2\,dX_2\\
		&\qquad+\int_\R p_n(t)\varphi(t)^2\,dt\int_\R L^{-2}|\eta'(X_2/L)|^2\,dX_2 .
	\end{align*}
	Here $p_n=0$ for $t\le-r_n$. For fixed $L$, the denominator converges uniformly on $[-L,L]$ to $|F_\infty'(t_+)|>0$. We keep $L$ fixed and let $n\to\infty$. The integrals with weight $p$ are over $I$ since $p=0$ outside $I$. Changing variables in the integrals involving $\eta$ and dividing by $L\|\eta\|_{L^2(\R)}^2$, we obtain, since $\varphi(t_+)=1$,
	\begin{equation*}
		\frac{15p(t_+)}{|F_\infty'(t_+)|}
		\le\int_I p(t)|\varphi'(t)|^2\,dt+
		\frac{\|\eta'\|_{L^2(\R)}^2}{L^2\|\eta\|_{L^2(\R)}^2}\int_I p(t)\varphi(t)^2\,dt .
	\end{equation*}
	As $L\to\infty$, the last term tends to zero. By \eqref{2-7}, the right-hand side is strictly smaller than the left-hand side for all sufficiently large $L$. This is a contradiction.
\end{proof}

\subsection{Volume, energy, and support}\label{sec2-3}

The volume estimate combines a Pohozaev identity in $\R^5$ with a lower bound on the potential at a point of large radius. For the latter estimate we first integrate the Newtonian kernel in the angular variable. In the meridional formulas below, write $x=(r,z)$ and $y=(r',z')$, and set
\begin{equation}\label{2-9}
	\rho(x,y):=\frac{(r-r')^2+(z-z')^2}{rr'}
\end{equation}
and
\begin{equation}\label{2-2}
	G(t):=\frac1{8\pi^2}\int_{\mathbb S^3} \bigl(t+2(1-\omega_1)\bigr)^{-3/2} \,d\omega .
\end{equation}
To obtain the meridional formula, evaluate the five-dimensional potential at $(re_1,z)$ and write the integration point as $(r'\omega,z')$ with $\omega\in\mathbb S^3$. The distance appearing in the Newton kernel is
\begin{equation*}
	|(re_1,z)-(r'\omega,z')|^2
	= (r-r')^2+(z-z')^2+2rr'(1-\omega_1)
	=rr'\bigl(\rho(x,y)+2(1-\omega_1)\bigr).
\end{equation*}
The volume element is $(r')^3\,dr'\,d\omega\,dz'$. Multiplying the factor $(rr')^{-3/2}$ from the kernel by $(r')^3$ yields $r^{-3/2}(r')^{3/2}$. Integrating in $\omega$, we have
\begin{equation}\label{2-3}
	v(r,z)=15r^{-3/2}\int_DG(\rho(x,y))(r')^{3/2}\,dr'\,dz'.
\end{equation}
The kernel is positive and decreasing, and
\begin{equation}\label{2-4}
	\begin{split}
		G(t)&=\frac1{4\pi}\ln\frac1t+c_0+O(t|\ln t|),\\
		G'(t)&=-\frac1{4\pi t}+O(1+|\ln t|)\qquad(t\downarrow0).
	\end{split}
\end{equation}
The expansions in \eqref{2-4} follow from Lemma~2.1 in \cite{GS15}, where the meridional kernel is $2\pi G$ in the notation of \eqref{2-2}.

\begin{proof}[Proof of Theorem~\ref{thm1-1}]
	We first derive the Pohozaev identity, writing $x,y\in\R^5$ in the five-dimensional integrals. By Lemmas~\ref{lem2-1} and \ref{lem2-5},
	\begin{equation}\label{2-38}
		\|v\|_{C^{1,\alpha}(\R^5)}\le C_{K,\alpha},\qquad
		D^{(5)}\subset\{|z|\le C_K\}.
	\end{equation}
	For $k>0$ the bounded core is separated from the axis, and for $k=0$ we have $F_+=(v-1)_+$. Thus $F_+\in W_c^{1,\infty}(\R^5)$. The differentiated Newton kernel has size $C|x-y|^{-4}$, which is integrable near $x=y$ in five dimensions. We may therefore differentiate the potential formula for $v$. Exchanging $x$ and $y$ and adding the resulting integrals, we obtain
	\begin{align*}
		2\int_{D^{(5)}}x\cdot\nabla v(x)\,dx &=-\frac{45}{8\pi^2}\iint_{D^{(5)}\times D^{(5)}} \frac{x\cdot(x-y)+y\cdot(y-x)}{|x-y|^5}\,dx\,dy\\
		&=-\frac{45}{8\pi^2}\iint_{D^{(5)}\times D^{(5)}} |x-y|^{-3}\,dx\,dy =-3\int_{D^{(5)}}v\,dx .
	\end{align*}
	On $D^{(5)}$ we have $v=F+1+k/r^2$ and $x\cdot\nabla(r^{-2})=-2r^{-2}$. Also $\nabla F_+=\boldsymbol 1_{D^{(5)}}\nabla F$ almost everywhere. Using these identities and $\operatorname{div}x=5$,
	\begin{equation*}
		\int_{D^{(5)}}x\cdot\nabla v\,dx
		=\int_{\R^5}x\cdot\nabla F_+\,dx-2k\int_{D^{(5)}}r^{-2}\,dx
		=-5\int_{D^{(5)}}F\,dx-2k\int_{D^{(5)}}r^{-2}\,dx.
	\end{equation*}
	Equating the two expressions for $2\int_{D^{(5)}}x\cdot\nabla v$,
	\begin{equation*}
		-3\int_{D^{(5)}}(F+1+k/r^2)\,dx=-10\int_{D^{(5)}}F\,dx-4k\int_{D^{(5)}}r^{-2}\,dx.
	\end{equation*}
	Thus
	\begin{equation}\label{2-43}
		7\int_{D^{(5)}}F\,dx+k\int_{D^{(5)}}r^{-2}\,dx=3|D^{(5)}|.
	\end{equation}

	To bound the volume in \eqref{2-43}, we show that $7F+k/r^2\le3/2$ on $D$ for large $r$. For $r^2\ge4K/3$, a point with $7F+k/r^2>3/2$ has $7F>3/4$. Let $x_0=(r_0,z_0)$ be such a point with $r_0\ge2$. Choose $\delta=\delta(K)<1$ with $7\delta\sup_{r\ge1}|\nabla F|<3/4$. Then the planar disk $B_\delta(x_0)$ is contained in $D$. For $y=(r',z')$ in this disk, $r_0/2\le r'\le3r_0/2$ and $\rho(x_0,y)\le C|x_0-y|^2/r_0^2$. By \eqref{2-3}--\eqref{2-4},
	\begin{equation*}
		C_K\ge v(r_0,z_0)
		\ge c\int_{B_\delta(x_0)}\ln\frac{r_0}{|x_0-y|}\,dy-C_K
		\ge c\pi\delta^2\ln r_0-C_K.
	\end{equation*}
	This bounds $r_0$ in terms of $K$, so $7F+k/r^2\le3/2$ on $D\cap\{r>R_0\}$ for some $R_0=R_0(K)$. On $D$ we have $7F+k/r^2\le 8\|v\|_{L^\infty(\R^5)}$. Also \eqref{2-38} implies $|D^{(5)}\cap\{r\le R_0\}|\le 2C_K|\mathbb S^3|\int_0^{R_0}r^3\,dr$. Splitting \eqref{2-43} at $R_0$, we obtain
	\begin{equation*}
		3|D^{(5)}|
		\le 8\|v\|_{L^\infty(\R^5)}|D^{(5)}\cap\{r\le R_0\}|
		+\frac{3}{2}|D^{(5)}\cap\{r>R_0\}|
		\le C_K+\frac{3}{2}|D^{(5)}|.
	\end{equation*}
	Moving the last term to the left bounds $|D^{(5)}|$ by $C_K$. Testing \eqref{1-1} with $v$ by density and using the volume bound and Sobolev's inequality, we obtain
	\begin{equation*}
		\int_{\R^5}|\nabla v|^2\,dx=15\int_{D^{(5)}}v\,dx\le C_K,
		\qquad \|v\|_{L^{10/3}(\R^5)}\le C_K.
	\end{equation*}
	To bound $r$ on $D$, take $(r_0,z_0)\in D$. Since $v(r_0,z_0)\ge1$, the gradient bound and axial symmetry imply $v(x',z)\ge1/2$ for $\bigl||x'|-r_0\bigr|<\delta$ and $|z-z_0|<\delta$ after decreasing $\delta(K)$. If $r_0\ge2\delta$,
	\begin{equation*}
		C_K\ge\int_{\R^5}v^{10/3}\,dx
		\ge 2^{-10/3}|\mathbb S^3|(2\delta)
		\int_{r_0-\delta}^{r_0+\delta}r^3\,dr
		\ge c\delta^2r_0^3.
	\end{equation*}
	Hence $r_0\le C_K$ for all $(r_0,z_0)\in D$.
\end{proof}

\subsection{The Hill limit and global continuation}\label{sec2-4}

For a sequence of solutions with bounded flux, Theorem~\ref{thm1-1} bounds the supports in a common ball. The limit is strictly decreasing in $z$ for $z>0$, which implies convergence of the characteristic functions. Testing the difference equation, we obtain strong convergence in the energy norm.

\begin{lemma}\label{lem2-6}
	Let $k_n\to k\ge0$ and let $v_n$ be centered finite-energy solutions of \eqref{1-1}. Write $D_n$ for the cross-section of $v_n$. A subsequence converges to a nonzero centered solution $v$ at $k$,
	\begin{equation*}
		\begin{aligned}
			v_n&\to v &&\text{in }D^{1,2}(\R^5)\cap C^{1,\beta}(\R^5) &&(0<\beta<1),\\
			\boldsymbol 1_{D_n^{(5)}}&\to\boldsymbol 1_{D^{(5)}} &&\text{in }L^q(\R^5) &&(1\le q<\infty).
		\end{aligned}
	\end{equation*}
\end{lemma}
\begin{proof}
	By Theorem~\ref{thm1-1}, $D_n^{(5)}\subset B_{R_*}(0)$ for a fixed $R_*$, and $v_n$ is bounded in $D^{1,2}(\R^5)\cap C^{1,\alpha}(\R^5)$ for each $0<\alpha<1$. Take a diagonal subsequence so that for every $0<\beta<1$,
	\begin{equation}\label{2-49}
		v_n\to v\quad\text{in }C^{1,\beta}_{\mathrm{loc}}(\R^5),\qquad
		\nabla v_n\rightharpoonup\nabla v\quad\text{in }L^2(\R^5).
	\end{equation}
	Take $x_n\in D_n^{(5)}$ and a subsequence with $x_n\to x_*$. Then $v(x_*)=\lim_n v_n(x_n)\ge1$. By the formula for the Newtonian potential, uniformly in $n$,
	\begin{equation*}
		0\le v_n(x)\le C_{R_*}|x|^{-3}\quad(|x|\ge2R_*).
	\end{equation*}
	The same bound holds for $v$. Together with \eqref{2-49} it implies $\|v_n-v\|_{L^\infty(\R^5)}\to0$.

	The limit $v$ is even in $z$. Since $\boldsymbol 1_{D_n^{(5)}}$ is nonincreasing in $z>0$,
	\begin{equation}\label{2-53}
		w_n=-\partial_zv_n\ge0,\qquad
		-\Delta w_n=-15\partial_z\boldsymbol 1_{D_n^{(5)}}\ge0
		\quad\text{in }\{z>0\}.
	\end{equation}
	The right-hand sides of \eqref{1-1} are uniformly bounded, so for each $1<p<\infty$ the sequence $(v_n)$ is bounded in $W^{2,p}_{\mathrm{loc}}(\R^5)$. In particular, $w=-v_z\in W^{1,p}_{\mathrm{loc}}(\R^5)$. By \eqref{2-49}, $w\ge0$ and $-\Delta w\ge0$ in $\{z>0\}$. If $w\equiv0$, evenness and decay force $v\equiv0$. The strong maximum principle implies
	\begin{equation}\label{2-54}
		-v_z=w>0\quad\text{throughout }\{z>0\},
	\end{equation}
	including the axis.

	For $r>0$, $F_n=v_n-1-k_n/r^2$ tends pointwise to $F=v-1-k/r^2$. By evenness and \eqref{2-54}, each section of $\{F=0\}$ at fixed $r>0$ has at most two points, so this set has measure zero by Fubini. Outside this set and the axis, the signs of $F_n$ agree with the sign of $F$ for large $n$. The common support bound then implies
	\begin{equation}\label{2-52}
		\boldsymbol 1_{D_n^{(5)}}\to\boldsymbol 1_{D^{(5)}}\quad\text{almost everywhere and in }L^q(\R^5)
		\quad(1\le q<\infty).
	\end{equation}
	Hence $-\Delta v=15\boldsymbol 1_{D^{(5)}}$. Subtracting the equations for $v_n$ and $v$, testing with $v_n-v$, and using Sobolev's inequality, we have
	\begin{align*}
		\|\nabla(v_n-v)\|_{L^2(\R^5)}^2 &=15\int_{\R^5}(\boldsymbol 1_{D_n^{(5)}}-\boldsymbol 1_{D^{(5)}})(v_n-v)\,dx\\
		&\le C\|\boldsymbol 1_{D_n^{(5)}}-\boldsymbol 1_{D^{(5)}}\|_{L^{10/7}(\R^5)} \|\nabla(v_n-v)\|_{L^2(\R^5)}.
	\end{align*}
	Thus $\|\nabla(v_n-v)\|_{L^2(\R^5)}\to0$. For $p>5/(1-\beta)$, apply \eqref{2-14} to $v_n-v$ on each unit ball. As in the proof of Lemma~\ref{lem2-1}, we use the interior estimate for $|x-y|<1/2$ and the supremum of the gradient for $|x-y|\ge1/2$ to obtain
	\begin{equation}\label{2-56}
		\|v_n-v\|_{C^{1,\beta}(\R^5)}
		\le C_p\left(\|v_n-v\|_{L^\infty(\R^5)}
		+\|\boldsymbol 1_{D_n^{(5)}}-\boldsymbol 1_{D^{(5)}}\|_{L^p(\R^5)}\right)\to0.
	\end{equation}
\end{proof}

The uniform energy bound enables us to apply the uniqueness theorem of Amick and Fraenkel for small flux.

\begin{corollary}\label{cor2-7}
	There is $k_0>0$ such that, for every $0\le k\le k_0$, \eqref{1-1} has exactly one nontrivial centered finite-energy solution. At $k=0$ it is Hill's vortex. For $0<k\le k_0$ it is Norbury's local solution.
\end{corollary}
\begin{proof}
	At $k=0$ the conclusion follows from the uniqueness of Hill's vortex \cite{AF}. For $k>0$ we compare the two energy norms. For $\psi=r^2v$, direct differentiation yields
	\begin{equation*}
		\frac{|\nabla\psi|^2}{r}
		=r^3|\nabla v|^2+2\partial_r(r^2v^2).
	\end{equation*}
	For smooth compactly supported $v$, the last term integrates to zero in $r$, since $r^2v^2$ vanishes at $r=0$ and at infinity. Integrating in the angular variables and using $|\mathbb S^3|=2\pi^2$, we obtain the second equality below. By density, the identity holds for $v\in D^{1,2}(\R^5)$,
	\begin{equation}\label{2-57}
		\int_0^\infty\int_\R\frac{|\nabla\psi|^2}{r}\,dz\,dr
		=\int_0^\infty\int_\R r^3|\nabla v|^2\,dz\,dr
		=\frac{1}{2\pi^2}\|\nabla v\|_{L^2(\R^5)}^2.
	\end{equation}
	The first integral is the square of the energy norm of the stream function defined on p.~209 of \cite{AF88}.

	Amick and Fraenkel use propagation speed $4/3$ and vortex strength $10$. Multiplication of both $\psi$ and $k$ by $2/3$ preserves the core and puts the equation in their normalization. By Theorem~\ref{thm1-1} with $K=1$ and \eqref{2-57}, the rescaled stream functions for $0<k\le1$ have a uniformly bounded energy norm in that normalization. We choose a bound larger than the norms of these solutions and of Hill's solution. Theorem~4.3 of \cite{AF88} yields the boundary estimate (1.12) of that paper for all such centered solutions at sufficiently small flux. Theorem~1.1 of the same paper then identifies them with Norbury's solution. The existence of that solution follows from Theorem~3.1 in \cite{Nor72}.
\end{proof}

We next use the energy bound to show that $k$ is unbounded on the Amick--Turner continuum.

\begin{corollary}\label{cor2-8}
	The set of fluxes on the Amick--Turner continuum from Hill's vortex is $[0,\infty)$.
\end{corollary}

\begin{proof}
	The continuum in \cite{AT88} is closed, connected and unbounded in the product of $[0,\infty)$ and the energy space of the stream function. We use the normalization of \eqref{1-1}. By the weak formulation (1.9) in \cite{AF88}, for each solution the function $v=\psi/r^2$ belongs to $D^{1,2}(\R^5)$ and satisfies
	\begin{equation*}
		-\Delta_{\R^5}v=15\boldsymbol 1_{D^{(5)}}
	\end{equation*}
	in $\R^5$. At $k=0$ we use Hill's solution after an axial translation \cite{AF}. For $k>0$, Lemma~2.1 of \cite{AF88} allows us to translate in $z$ so that $v$ is even and $v_z<0$ for $r,z>0$. Lemma~2.3 of the same paper applied to each such solution shows that $D$ is bounded.

	Since the right-hand side is bounded, $v\in W^{2,p}_{\mathrm{loc}}(\R^5)$ for $1<p<\infty$. The Newtonian potential of $15\boldsymbol 1_{D^{(5)}}$ belongs to $D^{1,2}(\R^5)$ and decays uniformly at infinity, and hence is equal to $v$ by uniqueness in $D^{1,2}(\R^5)$. In $\{z>0\}$ the characteristic function of $D^{(5)}$ is nonincreasing in $z$, and by continuity $-v_z\ge0$ also on the axis. Thus
	\begin{equation*}
		-\Delta(-v_z)=-15\partial_z\boldsymbol 1_{D^{(5)}}\ge0
		\quad\text{in }\{z>0\}.
	\end{equation*}
	Since $-v_z$ is already positive for $r,z>0$, the strong maximum principle implies $-v_z>0$ throughout $\{z>0\}$.

	Theorem~\ref{thm1-1} now applies. Axial translations preserve the flux and the energy norm. If $k\le K$ on the continuum, \eqref{2-57} bounds the energy norm of $\psi$ uniformly on the continuum. This contradicts its unboundedness. Since the continuum is connected and contains Hill's vortex at $k=0$, its fluxes form an unbounded interval in $[0,\infty)$ containing $0$. Hence this interval is $[0,\infty)$.
\end{proof}

\begin{remark}\label{rem2-1}
	It remains open whether two points of the continuum can have the same intermediate flux.
\end{remark}

\section{The radial location of the cross-section}\label{sec3}

Cao et al.~\cite{CWWZ} and Cao and Zhan \cite{CZ} determine the limiting radius in the whole space by comparing energies at prescribed circulation. At fixed flux the reference radius is $R=\sqrt{k/3}$, where the derivative of $B$ in \eqref{3-1} vanishes. We show that the radial coordinate of every point of $D$, divided by $R$, tends to one. We bound the smallest radius of $D$ using the sign of the radial derivative of the kernel for points on opposite sides of $r=T$. To bound the largest radius, we first estimate the mass on $\{r>(1+\delta)R\}$ and then apply the integral identity on that side.

For a centered solution $(k,v)$ with $R\ge1$, Lemma~\ref{lem2-1} implies
\begin{equation}\label{2-60}
	\|v\|_{L^\infty(\R^5)}\le CR,\qquad \|\nabla v\|_{L^\infty(\R^5)}\le C\sqrt R.
\end{equation}
Since $F_+\le v$, Lemma~\ref{lem2-5} and \eqref{2-60} imply
\begin{equation}\label{2-62}
	\sup_D|z|\le C\sqrt R=o(R).
\end{equation}

\subsection{Integral identities for the radial estimates}\label{sec3-1}

We bound the smallest and largest values of $r$ on $D$ by integrating on either side of $r=T$. Multiplying $v$ by $r^{3/2}$ puts the potential formula in a symmetric form. Set
\begin{equation}\label{3-1}
	V:=r^{3/2}v,\quad B(r):=r^{3/2}+kr^{-1/2},\quad V-B=r^{3/2}F .
\end{equation}
Thus $V=B$ on the free boundary and $(V-B)_+=r^{3/2}F_+$. Differentiating $B$, we have
\begin{equation*}
	B'(r)=\frac32r^{-3/2}(r^2-R^2).
\end{equation*}
Thus $B'$ changes sign at the reference radius $R$. With $d\mu=15r^{3/2}\boldsymbol 1_D\,dr\,dz$, \eqref{2-3} reads
\begin{equation}\label{3-2}
	V(x)=\int_{D}G(\rho(x,y))\,d\mu(y).
\end{equation}
For each fixed solution, $D\subset\{v>1\}$ is bounded by decay. The inequality $v>k/r^2$ on $D$ implies
\begin{equation}\label{3-4}
	\overline D\subset\{r\ge\sqrt{k/\|v\|_{L^\infty(\R^5)}}\}.
\end{equation}
In particular,
\begin{equation}\label{3-5}
	r^{3/2}F_+\in W^{1,p}_c((0,\infty)\times\R),\qquad
	\nabla(r^{3/2}F_+)=\boldsymbol 1_D\nabla(r^{3/2}F).
\end{equation}
The radial identities also use the global derivative bounds
\begin{equation}\label{2-5}
	0<-tG'(t)\le C,\qquad
	-tG'(t)\le\frac{3}{2}G(t)\qquad(t>0).
\end{equation}
The first bound follows from Remark~2.2 in \cite{GS15}. By \eqref{2-2},
\begin{equation}\label{3-8}
	-tG'(t)=\frac3{16\pi^2}\int_{\mathbb S^3}
	\frac{t}{\bigl(t+2(1-\omega_1)\bigr)^{5/2}}\,d\omega
	\le\frac32G(t),
\end{equation}
which proves the second bound in \eqref{2-5}.

Since $r$ is bounded above and away from zero on $D$, the identity
\begin{equation*}
	\partial_r\rho=\frac{2(r-r')}{rr'}-\frac{\rho}{r}
\end{equation*}
implies $|\partial_r\rho|\le C_D(|x-y|+|x-y|^2)$. Since $|G'(\rho)|\le C/\rho$ and $\rho=|x-y|^2/(rr')$,
\begin{equation*}
	|r^{-3/2}\partial_rG(\rho(x,y))|
	\le C_D(1+|x-y|^{-1})\in L^1(D\times D).
\end{equation*}
The density of $d\mu$ is bounded on $D$. We differentiate \eqref{3-2} in $r$ and average the resulting double integral with the one obtained by exchanging $x=(r,z)$ and $y=(r',z')$. Since $\partial_rG(\rho)$ is integrable on $D\times D$,
\begin{equation}\label{3-7}
	\int_D r^{-3/2}V_r\,d\mu
	=\frac12\iint_{D\times D}G'(\rho(x,y))
	(r^{-3/2}\partial_r+(r')^{-3/2}\partial_{r'})\rho
	\,d\mu(x)\,d\mu(y).
\end{equation}
The derivatives of $\rho$ are
\begin{equation}\label{3-10}
	\begin{aligned}
		\partial_r\rho&=\frac{r^2-(r')^2-(z-z')^2}{r^2r'},& \partial_z\rho&=\frac{2(z-z')}{rr'},\\
		\partial_{r'}\rho&=\frac{(r')^2-r^2-(z-z')^2}{r(r')^2},& \partial_{z'}\rho&=-\frac{2(z-z')}{rr'}.
	\end{aligned}
\end{equation}
Consequently, for $x\neq y$,
\begin{equation}\label{3-9}
	\begin{split}
		&(r^{-3/2}\partial_r+(r')^{-3/2}\partial_{r'})G(\rho)\\
		&\quad=[-\rho G'(\rho)]\left[ r^{-5/2}+(r')^{-5/2} +3\frac{(r-r')^2}{|x-y|^2} \int_0^1(r'+\tau(r-r'))^{-5/2}\,d\tau \right].
	\end{split}
\end{equation}
Indeed, differentiating $\rho=|x-y|^2/(rr')$,
\begin{equation*}
	(r^{-3/2}\partial_r+(r')^{-3/2}\partial_{r'})\rho
	=\frac{2(r-r')(r^{-3/2}-(r')^{-3/2})}{rr'}
	-(r^{-5/2}+(r')^{-5/2})\rho,
\end{equation*}
and $r^{-3/2}-(r')^{-3/2} =-(3/2)(r-r')\int_0^1(r'+\tau(r-r'))^{-5/2}\,d\tau$. The integral in \eqref{3-9} is at most $\max(r^{-5/2},(r')^{-5/2})$. Since $(r-r')^2\le|x-y|^2$ and $0<-\rho G'(\rho)\le C$ by \eqref{2-5}, the expression in \eqref{3-9} lies between $0$ and $C(r^{-5/2}+(r')^{-5/2})$.

We divide $D$ by the line $r=T$ and write
\begin{align*}
	D^-(T)&:=D\cap\{r<T\},& D^+(T)&:=D\cap\{r>T\},\\
	S^-(T)&:=|D^-(T)|,& S^+(T)&:=|D^+(T)|.
\end{align*}
Writing $H(r)$ as in the proof of Lemma~\ref{lem2-2}, we have
\begin{equation*}
	D\cap(\{r\}\times\R)=\{r\}\times(-H(r),H(r)), \qquad H(r)=0\quad\text{for an empty section}.
\end{equation*}
We integrate $\partial_r(r^{3/2}F_+)$ over $0<r<T$ and over $T<r<\infty$ to express the integral of $F_+(T,z)$ in two ways.

\begin{lemma}\label{lem3-1}
	For almost every $T>0$,
	\begin{equation}\label{3-12}
		T^{3/2}\int_{\R}F_+(T,z)\,dz
		=\int_{D^-(T)}(V_r-B')\,dr\,dz
		=\int_{D^+(T)}(B'-V_r)\,dr\,dz .
	\end{equation}
	Moreover, $S^-$ and $S^+$ are absolutely continuous and
	\begin{equation}\label{3-13}
		(S^-)'(T)=2H(T),\qquad -(S^+)'(T)=2H(T) \quad\text{for almost everywhere }T.
	\end{equation}
\end{lemma}
\begin{proof}
	For almost every $z$, the function $(V-B)_+(r,z)$ is absolutely continuous in $r$ and vanishes near $r=0$ and at infinity by \eqref{3-5}. Integrating its derivative on the two sides of $T$,
	\begin{align*}
		(V-B)_+(T,z) &=\int_0^T\boldsymbol 1_D(r,z)\bigl(V_r(r,z)-B'(r)\bigr)\,dr\\
		&=\int_T^\infty\boldsymbol 1_D(r,z)\bigl(B'(r)-V_r(r,z)\bigr)\,dr.
	\end{align*}
	Integrating in $z$ and using $(V-B)_+(T,z)=T^{3/2}F_+(T,z)$ proves \eqref{3-12}. For \eqref{3-13}, differentiate
	\begin{equation*}
		S^-(T)=2\int_0^T H(r)\,dr,\qquad S^+(T)=2\int_T^\infty H(r)\,dr.
	\end{equation*}
\end{proof}

For $T>\sup_D r$, the left-hand side of \eqref{3-12} is zero and $D^-(T)=D$. Thus $\int_D B'\,dr\,dz=\int_D V_r\,dr\,dz$. By \eqref{3-7} and $d\mu=15r^{3/2}\,dr\,dz$ on $D$, this becomes
\begin{equation}\label{3-14}
	\begin{split}
		\int_D B'(r)\,dx &=\frac1{15}\int_D r^{-3/2}V_r\,d\mu\\
		&=\frac{15}{2}\iint_{D\times D}r^{3/2}(r')^{3/2} (r^{-3/2}\partial_r+(r')^{-3/2}\partial_{r'})G(\rho) \,dx\,dy,
	\end{split}
\end{equation}
where $dx=dr\,dz$ and $dy=dr'\,dz'$.

For a nonempty section, evenness and $F(r,H(r))=0$ imply, for $|z|\le H(r)$,
\begin{equation}\label{3-16}
	\begin{split}
		F(r,z) &=-\int_{|z|}^{H(r)}F_z(r,t)\,dt\\
		&\le\|F_z(r,\cdot)\|_{L^\infty(-H(r),H(r))} (H(r)-|z|),\\
		\int_\R F_+(r,z)\,dz &\le\|F_z(r,\cdot)\|_{L^\infty(-H(r),H(r))}H(r)^2.
	\end{split}
\end{equation}
The bound on $\int_\R F_+(r,z)\,dz$ also holds at $H(r)=0$. If the left-hand side of \eqref{3-12} is bounded below by a multiple of $S^\pm(T)$, then \eqref{3-16} implies $H(T)^2\ge cS^\pm(T)$. We integrate this inequality using \eqref{3-13}.

\begin{lemma}\label{lem3-2}
	Let $S\ge0$ be absolutely continuous on $[T_{0},T_{1}]$ and positive in its interior, and let $H\ge0$ with $H^2\ge cS$ almost everywhere for some $c>0$. Assume that either $S'=2H$ holds almost everywhere on the interval or $S'=-2H$ holds almost everywhere on the interval. Then
	\begin{equation}\label{3-18}
		\begin{split}
			|\sqrt{S(T_{1})}-\sqrt{S(T_{0})}| &\ge\sqrt{c}(T_{1}-T_{0}),\\
			|S(T_{1})-S(T_{0})| &\ge c(T_{1}-T_{0})^2.
		\end{split}
	\end{equation}
\end{lemma}
\begin{proof}
	On each closed interval inside $(T_0,T_1)$, the function $\sqrt S$ is absolutely continuous and monotone with
	\begin{equation*}
		\bigl|(\sqrt S)'(T)\bigr|=\frac{H(T)}{\sqrt{S(T)}}\ge\sqrt c
		\qquad\text{for almost every }T.
	\end{equation*}
	Integrating on closed subintervals of $(T_0,T_1)$ and letting the endpoints tend to $T_0$ and $T_1$, we obtain the first inequality in \eqref{3-18}. The second follows from
	\begin{equation*}
		|S(T_1)-S(T_0)|=
		|\sqrt{S(T_1)}-\sqrt{S(T_0)}|
		\bigl(\sqrt{S(T_1)}+\sqrt{S(T_0)}\bigr).
	\end{equation*}
\end{proof}

For $cR\le T\le CR$, \eqref{3-16} and the gradient bound \eqref{2-60} imply
\begin{equation}\label{3-23}
	T^{3/2}\int_\R F_+(T,z)\,dz
	\le CT^{3/2}\sqrt R\,H(T)^2
	\le CR^2H(T)^2.
\end{equation}
Let $cR\le T_0<T_1\le CR$ with $T_1-T_0\ge cR$. Take $S=S^-$ or $S=S^+$ and suppose $S$ is positive on $(T_0,T_1)$ and
\begin{equation*}
	T^{3/2}\int_\R F_+(T,z)\,dz\ge c\sqrt R\,S(T)
	\quad\text{for almost every }T\in(T_0,T_1).
\end{equation*}
Then $H(T)^2\ge cR^{-3/2}S(T)$. By \eqref{3-13} and Lemma~\ref{lem3-2},
\begin{equation}\label{3-24}
	\begin{split}
		|S(T_1)-S(T_0)| &\ge cR^{-3/2}(T_1-T_0)^2\ge cR^{1/2},\\
		\mu(D\cap\{T_0<r<T_1\}) &\ge15T_0^{3/2}|S(T_1)-S(T_0)|\ge cR^2.
	\end{split}
\end{equation}
\subsection{The inner radial estimate}\label{sec3-2}

For $r<T<R$ we have $-B'>0$. To bound the right-hand side of \eqref{3-12} from below, we show that $\int_{D^-(T)}V_r\,dr\,dz$ is nonnegative.

\begin{lemma}\label{lem3-3}
	For almost every $T<R$,
	\begin{equation}\label{3-17}
		T^{3/2}\int_\R F_+(T,z)\,dz \ge\int_{D^-(T)}-B'(r)\,dr\,dz.
	\end{equation}
\end{lemma}
\begin{proof}
	Write $dx=dr\,dz$. Since $D\cap\{r=T\}$ has area zero, we split the integral in \eqref{3-2} over $D^-(T)$ and $D^+(T)$. Using $15\,\,dx=r^{-3/2}\,d\mu(x)$ on $D$, we have
	\begin{align*}
		15\int_{D^-(T)}V_r\,dx &=\iint_{D^-(T)\times D^-(T)} r^{-3/2}\partial_rG(\rho)\,d\mu(x)\,d\mu(y)\\
		&\quad+\iint_{D^-(T)\times D^+(T)} r^{-3/2}\partial_rG(\rho)\,d\mu(x)\,d\mu(y).
	\end{align*}
	By exchanging $x,y$ and averaging, we write the first integral as one half of the integral of \eqref{3-9} over $D^-(T)\times D^-(T)$. It is therefore nonnegative. In the second integral $r<r'$, so $\partial_r\rho<0$ by \eqref{3-10}. Since $G'<0$, the integrand is nonnegative. Inserting $\int_{D^-(T)}V_r\,dx\ge0$ into \eqref{3-12} proves \eqref{3-17}.
\end{proof}

For the lower bound on $r$, we will take $R/2<T<R$. To bound the mass of $D$ in $\{R/2<r<2R\}$, we evaluate $V$ at a boundary point in that region.

\begin{lemma}\label{lem3-4}
	For all sufficiently large $k$,
	\begin{equation}\label{3-19}
		\mu\bigl(D\cap\{R/2<r<2R\}\bigr)\le CR^{3/2}.
	\end{equation}
\end{lemma}

\begin{proof}
	It suffices to consider the case $D\cap\{R/2<r<2R\}\ne\varnothing$. Take $x_0=(r_0,H(r_0))$ with $H(r_0)>0$ and $R/2<r_0<2R$. For $y=(r',z')\in D$ with $R/2<r'<2R$, we have
	\begin{align*}
		\rho(x_0,y)&=\frac{(r_0-r')^2+(H(r_0)-z')^2}{r_0r'}\\
		&\le\frac{(3R/2)^2+4\sup_D|z|^2}{(R/2)^2}\le C,
	\end{align*}
	where we have used $\sup_D|z|\le R$ for large $k$, by \eqref{2-62}. In view of \eqref{3-2} and $V=B$ on the boundary of $D$, we have
	\begin{align*}
		G(C)\mu\bigl(D\cap\{R/2<r<2R\}\bigr)&\le\int_D G(\rho(x_0,y))\,d\mu(y)\\
		&=V(x_0)=B(r_0)\\
		&=r_0^{3/2}+3R^2r_0^{-1/2}\\
		&\le(2R)^{3/2}+3R^2(R/2)^{-1/2}\le CR^{3/2}.
	\end{align*}
	In the first inequality we used $G(\rho(x_0,y))\ge G(C)>0$ for $R/2<r'<2R$.
\end{proof}

If $D$ contains a point with $r<(1-2\delta)R$, we apply Lemma~\ref{lem3-3} for $(1-2\delta)R<T<(1-\delta)R$. The resulting lower bound for the mass of $D$ between these two radii contradicts Lemma~\ref{lem3-4}.

\begin{corollary}\label{cor3-5}
	For every $0<\delta<1/4$ and all sufficiently large $k$,
	\begin{equation*}
		D\subset\{r\ge(1-2\delta)R\}.
	\end{equation*}
\end{corollary}
\begin{proof}
	Suppose $D\cap\{r<(1-2\delta)R\}\neq\varnothing$. Set
	\begin{equation*}
		T_{0}=(1-2\delta)R,\qquad T_{1}=(1-\delta)R.
	\end{equation*}
	For $T_0<T<T_1$ we have $S^-(T)>0$. Since
	\begin{equation*}
		-B'(r)=\frac32\frac{R^2-r^2}{r^{3/2}}\ge c_\delta\sqrt R\qquad(0<r\le T_{1}),
	\end{equation*}
	Lemma~\ref{lem3-3} implies
	\begin{equation*}
		T^{3/2}\int_\R F_+(T,z)\,dz\ge c_\delta\sqrt R\,S^-(T).
	\end{equation*}
	Since $T_1-T_0=\delta R$, \eqref{3-24} bounds the mass of $D\cap\{T_0<r<T_1\}$ below by $c_\delta R^2$. Since $R/2<T_0<T_1<2R$, Lemma~\ref{lem3-4} implies
	\begin{equation*}
		\mu(D\cap\{T_0<r<T_1\})
		\le\mu(D\cap\{R/2<r<2R\})\le CR^{3/2},
	\end{equation*}
	which contradicts \eqref{3-24} as $R\to\infty$.
\end{proof}

\subsection{Mass away from the reference radius}\label{sec3-4}

On $r>(1+\delta)R$, the function $rB'$ is bounded below by a positive multiple of $B$. To estimate the $\mu$-measure of this part of $D$, we first bound $\int_D rB'\,d\mu$ by the Pohozaev identity and the formula for $V_z$.

\begin{lemma}\label{lem3-6}
	For a centered finite-energy solution of \eqref{1-1} with $k>0$,
	\begin{equation}\label{3-21}
		\begin{split}
			\int_DrB'\,d\mu
			&=\frac72\int_D r^{3/2}F\,d\mu\\
			&=\frac72\iint_{D\times D}[-G'(\rho)]\frac{(z-z')^2}{rr'}\,d\mu(x)\,d\mu(y)\ge0 .
		\end{split}
	\end{equation}
	For all sufficiently large $k$,
	\begin{equation}\label{3-22}
		\int_D rB'\,d\mu\le\frac{C}{\ln R}\int_D V\,d\mu .
	\end{equation}
\end{lemma}

\begin{proof}
	Dividing the Pohozaev identity \eqref{2-43} by $|\mathbb S^3|$, we obtain
	\begin{equation*}
		7\int_DFr^3\,dr\,dz+k\int_Dr\,dr\,dz
		=3\int_Dr^3\,dr\,dz .
	\end{equation*}
	Using $d\mu=15r^{3/2}\,dr\,dz$ on $D$,
	\begin{align*}
		\int_DrB'\,d\mu
		&=\frac{15}{2}\left(3\int_Dr^3\,dr\,dz-k\int_Dr\,dr\,dz\right)\\
		&=\frac{105}{2}\int_DFr^3\,dr\,dz
		=\frac72\int_Dr^{3/2}F\,d\mu .
	\end{align*}
	Since $F_+$ has compact support and $(F_+)_z=\boldsymbol 1_DF_z$, integration by parts in $z$ yields
	\begin{align*}
		\int_Dr^{3/2}F\,d\mu
		&=15\int_{(0,\infty)\times\R} r^3F_+\,dr\,dz\\
		&=-15\int_D zr^3F_z\,dr\,dz
		=-\int_D zV_z\,d\mu .
	\end{align*}
	Since $r,r'$ are bounded above and away from zero by \eqref{3-4}, \eqref{2-5} and \eqref{3-10} imply $|\partial_zG(\rho)|\le C_D|x-y|^{-1}$. The differentiated kernel is locally integrable, and its product with $z$ is integrable on $D\times D$ with respect to $d\mu(x)\,d\mu(y)$. We may therefore differentiate \eqref{3-2} in the weak sense and exchange $x$ and $y$ in the resulting double integral. Using \eqref{3-10},
	\begin{align*}
		-\int_D zV_z\,d\mu
		&=-\frac12\iint_{D\times D}G'(\rho)(z\partial_z+z'\partial_{z'})\rho\,d\mu(x)\,d\mu(y)\\
		&=\iint_{D\times D}[-G'(\rho)]\frac{(z-z')^2}{rr'}\,d\mu(x)\,d\mu(y).
	\end{align*}
	This proves \eqref{3-21}.

	For large $k$, $D\subset\{r\ge3R/4\}$ by Corollary~\ref{cor3-5}. By the height bound \eqref{2-62},
	\begin{equation*}
		\frac{(z-z')^2}{rr'}\le\min\{\rho,C/R\}\qquad(x,y\in D).
	\end{equation*}
	For $0<\rho\le R^{-1/2}$, \eqref{2-4} and \eqref{2-5} imply
	\begin{equation*}
		G(\rho)\ge c\ln R,\qquad
		[-G'(\rho)]\frac{(z-z')^2}{rr'}\le-\rho G'(\rho)\le C.
	\end{equation*}
	For $\rho>R^{-1/2}$, the second bound in \eqref{2-5} implies
	\begin{equation*}
		[-G'(\rho)]\frac{(z-z')^2}{rr'}
		\le \frac{C}{R\rho}G(\rho)\le \frac{C}{\sqrt R}G(\rho).
	\end{equation*}
	Combining the two estimates and using $R^{-1/2}\ln R\to0$, we have
	\begin{equation}\label{3-34}
		[-G'(\rho)]\frac{(z-z')^2}{rr'}\le\frac{C}{\ln R}G(\rho).
	\end{equation}
	Integrating this bound in \eqref{3-21} and using \eqref{3-2}, we obtain \eqref{3-22}.
\end{proof}

To use \eqref{3-22}, we estimate $\int_D V\,d\mu$ by the integral of $rB'$ and the integral over $D\cap\{R/2<r<2R\}$.

\begin{lemma}\label{lem3-7}
	For centered finite-energy solutions with $k\to\infty$, one has for every $\delta>0$
	\begin{equation}\label{3-37}
		\frac1{R^{3/2}}\mu\left(D\cap
		\left\{\left|\frac rR-1\right|>\delta\right\}\right)\to0.
	\end{equation}
\end{lemma}

\begin{proof}
	By Corollary~\ref{cor3-5}, $D\subset\{r\ge3R/4\}$ for all sufficiently large $k$. Since $B\le CR^{3/2}$ on $[3R/4,2R]$, Lemma~\ref{lem3-4} implies
	\begin{equation}\label{3-39}
		\int_{D\cap\{r\le2R\}}B\,d\mu
		\le CR^{3/2}\mu(D\cap\{R/2<r<2R\})\le CR^3.
	\end{equation}
	From the definition of $B$ we have
	\begin{equation*}
		\frac{rB'(r)}{B(r)}=\frac{3(r^2-R^2)}{2(r^2+3R^2)},
		\qquad
		\frac{d}{dr}\frac{rB'(r)}{B(r)}
		=\frac{12rR^2}{(r^2+3R^2)^2}>0 .
	\end{equation*}
	In particular, the ratio $rB'/B$ is increasing. Thus $B\le CrB'$ for $r>2R$, and for $r<R$ we have $rB'<0$ and $|rB'|\le B/2$. By \eqref{3-39}, we obtain
	\begin{align*}
		\int_D B\,d\mu
		&\le CR^3+C\int_{D\cap\{r>2R\}}rB'\,d\mu\\
		&\le CR^3+C\int_D rB'\,d\mu+
		C\int_{D\cap\{r<R\}}|rB'|\,d\mu\\
		&\le CR^3+C\int_D rB'\,d\mu .
	\end{align*}
	By $V=B+r^{3/2}F$ and \eqref{3-21},
	\begin{equation*}
		\int_D V\,d\mu=\int_D B\,d\mu+\frac27\int_D rB'\,d\mu .
	\end{equation*}
	Consequently, \eqref{3-22} implies
	\begin{equation*}
		\int_D rB'\,d\mu\le\frac{C}{\ln R}\left(R^3+\int_D rB'\,d\mu\right).
	\end{equation*}
	For large $R$, $C/\ln R\le1/2$. Moving $(C/\ln R)\int_D rB'\,d\mu$ to the left, we obtain $\int_DrB'\,d\mu\le CR^3/\ln R$.

	To estimate the mass away from $R$, observe that on $D\cap\{r<R\}$,
	\begin{equation*}
		\frac{|rB'|}{B}\le\frac12\left(1-\frac{r^2}{R^2}\right)\le 1-\frac rR=o(1)
	\end{equation*}
	uniformly by Corollary~\ref{cor3-5}. By \eqref{3-39}, $\int_{D\cap\{r<R\}}|rB'|\,d\mu=o(R^3)$. For fixed $\delta>0$, the ratio $rB'/B$ is bounded below by $c_\delta>0$ when $r>(1+\delta)R$. Therefore
	\begin{equation}\label{3-38}
		c_\delta\int_{D\cap\{r>(1+\delta)R\}} B\,d\mu
		\le\int_D rB'\,d\mu+\int_{D\cap\{r<R\}}|rB'|\,d\mu=o(R^3).
	\end{equation}
	Since $B\ge B(R)=4R^{3/2}$, the mass of $D\cap\{r>(1+\delta)R\}$ is $o(R^{3/2})$. By Corollary~\ref{cor3-5}, there are no points of $D$ with $r<(1-\delta)R$ for all sufficiently large $k$. This proves \eqref{3-37}.
\end{proof}

\subsection{The outer radial estimate}\label{sec3-5}

To bound the outer radius, we estimate $\int_{D^+(T)}V_r\,dr\,dz$ from above in \eqref{3-12}. The next two lemmas bound the positive contribution of this integral using the mass for $r'>(1+\delta/4)R$. By Corollary~\ref{cor3-5}, all points of $D$ have $r\ge R/2$ for large $k$.

\begin{lemma}\label{lem3-8}
	Let $E,E'\subset D\cap\{r\ge R/2\}$ be measurable. Suppose $r>r'$ whenever $(r,z)\in E$ and $(r',z')\in E'$. Then
	\begin{align}
		0\le\iint_{E\times E}r^{-3/2}\partial_rG(\rho(x,y))\,d\mu(x)\,d\mu(y)
		&\le CR^{-3/2}\mu(E)\int_E\frac{d\mu}{r},\label{3-25}\\
		\left[\iint_{E\times E'}r^{-3/2}\partial_rG(\rho(x,y))\,d\mu(x)\,d\mu(y)\right]_+
		&\le CR^{-3/2}\mu(E')\int_E\frac{d\mu}{r},\label{3-26}
	\end{align}
	where $C$ is universal.
\end{lemma}

\begin{proof}
	By symmetry and \eqref{3-9},
	\begin{equation*}
		0\le\iint_{E\times E}r^{-3/2}\partial_rG(\rho)\,d\mu(x)\,d\mu(y)
		\le C\mu(E)\int_E r^{-5/2}\,d\mu.
	\end{equation*}
	For $r>r'$, \eqref{3-10} implies
	\begin{equation*}
		[r^{-3/2}\partial_r G(\rho)]_+
		\le [-G'(\rho)]\frac{(z-z')^2}{r^{7/2}r'}
		\le [-\rho G'(\rho)]r^{-5/2}\le Cr^{-5/2}.
	\end{equation*}
	Integrate on $E\times E'$. In both estimates, $r^{-5/2}\le2^{3/2}R^{-3/2}/r$.
\end{proof}

For $T>(1+\delta/2)R$, we divide the integral in $y$ at $r'=(1+\delta/4)R$. The mass on $\{r'>(1+\delta/4)R\}$ is small by Lemma~\ref{lem3-7}. For $r'\le(1+\delta/4)R$ and $r>T$, we have $r-r'\ge\delta R/4$. The height bound then implies that the radial derivative of the kernel is nonpositive.

\begin{lemma}\label{lem3-9}
	Fix $\delta>0$ and suppose $\sup_D|z|=o(R)$. Assume also that
	\begin{equation}\label{3-40}
		\mu\bigl(D\cap\{r>(1+\delta/4)R\}\bigr)=o(R^{3/2}).
	\end{equation}
	For all sufficiently large $k$ and almost every $T\in((1+\delta/2)R,(1+2\delta)R)$,
	\begin{equation}\label{3-41}
		T^{3/2}\int_\R F_+(T,z)\,dz
		\ge c_\delta\sqrt R\,S^+(T).
	\end{equation}
\end{lemma}

\begin{proof}
	Set $T_*=(1+\delta/4)R$. Using $15\,dr\,dz=r^{-3/2}\,d\mu$ on the right-hand side of \eqref{3-12}, we obtain
	\begin{equation}\label{3-42}
		15T^{3/2}\int_\R F_+(T,z)\,dz
		=\int_{D^+(T)}r^{-3/2}B'\,d\mu
		-\int_{D^+(T)}r^{-3/2}V_r\,d\mu.
	\end{equation}
	For $x=(r,z)\in D^+(T)$ and $y=(r',z')\in D^-(T_*)$,
	\begin{equation*}
		r-r'\ge\delta R/4,
		\qquad r^2-(r')^2\ge c_\delta R^2,
		\qquad (z-z')^2=o(R^2).
	\end{equation*}
	Thus \eqref{3-10} implies, for large $k$,
	\begin{equation}\label{3-43}
		\begin{split}
			r^{-3/2}\partial_rG(\rho(x,y))&\le0,\\
			\iint_{D^+(T)\times D^-(T_*)}r^{-3/2}\partial_rG(\rho(x,y))\,d\mu(x)\,d\mu(y)&\le0.
		\end{split}
	\end{equation}
	Since $D\cap\{r=T\}$ has area zero, we can write
	\begin{equation*}
		D^+(T_*)=D^+(T)\cup(D\cap\{T_*<r<T\}).
	\end{equation*}
	Apply \eqref{3-25} to the integral with both points in $D^+(T)$ and \eqref{3-26} to the integral for $T_*<r'<T$. We obtain
	\begin{align*}
		&\left[\iint_{D^+(T)\times D^+(T_*)}r^{-3/2}\partial_rG(\rho) \,d\mu(x)\,d\mu(y)\right]_+\\
		&\quad\le CR^{-3/2} \bigl[\mu(D^+(T))+\mu(D\cap\{T_*<r<T\})\bigr] \int_{D^+(T)}\frac{d\mu}{r}\\
		&\quad=CR^{-3/2}\mu(D^+(T_*))\int_{D^+(T)}\frac{d\mu}{r} =o(1)\int_{D^+(T)}\frac{d\mu}{r}.
	\end{align*}
	For $r>T$,
	\begin{equation*}
		r^{-3/2}B'(r)=\frac3{2r}\left(1-\frac{R^2}{r^2}\right)
		\ge\frac{c_\delta}{r}.
	\end{equation*}
	By \eqref{3-40}, for large $k$ the coefficient $CR^{-3/2}\mu(D^+(T_*))$ is at most $c_\delta/2$. By \eqref{3-43}, the integral with $y\in D^-(T_*)$ is nonpositive. We now use these two bounds on the right-hand side of \eqref{3-42}. Since $r>T>R$ on $D^+(T)$, we have
	\begin{align*}
		15T^{3/2}\int_{\R} F_+(T,z)\,dz
		&\ge \frac{c_\delta}{2}\int_{D^+(T)}\frac{d\mu}{r}\\
		&=\frac{15c_\delta}{2}\int_{D^+(T)}\sqrt r\,dr\,dz\\
		&\ge \frac{15c_\delta}{2}\sqrt R\,S^+(T).
	\end{align*}
\end{proof}

The lower bound \eqref{3-41} is uniform for $T$ in an interval of length $3\delta R/2$. If $S^+(T)$ is positive throughout that interval, the lower bound in \eqref{3-24} contradicts \eqref{3-40}.

\begin{proposition}\label{pro3-10}
	Under the hypotheses of Lemma~\ref{lem3-9}, for all sufficiently large $k$,
	\begin{equation*}
		D\subset\{r\le(1+2\delta)R\}.
	\end{equation*}
\end{proposition}
\begin{proof}
	Suppose $D\cap\{r>(1+2\delta)R\}\neq\varnothing$. Put
	\begin{equation*}
		T_0=(1+\delta/2)R,\qquad T_1=(1+2\delta)R.
	\end{equation*}
	For $T_0<T<T_1$, $S^+(T)>0$. Since $T_1-T_0=3\delta R/2$, we can apply \eqref{3-24} to $S^+$ using \eqref{3-41}. Since $T_0>(1+\delta/4)R$, \eqref{3-40} implies
	\begin{equation*}
		c_\delta R^2
		\le\mu(D\cap\{T_0<r<T_1\})
		\le\mu(D\cap\{r>(1+\delta/4)R\})
		=o(R^{3/2}),
	\end{equation*}
	a contradiction as $R\to\infty$.
\end{proof}

\begin{proposition}\label{pro3-11}
	For centered finite-energy solutions $(k,v)$, as $k\to\infty$,
	\begin{equation}\label{3-44}
		\sup_{(r,z)\in D}\left|\frac r{R}-1\right|\to0, \qquad \sup_{D}\frac{|z|}{R}\to0.
	\end{equation}
\end{proposition}

\begin{proof}
	By \eqref{2-62}, $\sup_D|z|/R\to0$. Fix $0<\delta<1/4$. Lemma~\ref{lem3-7} implies \eqref{3-40}. Combining the lower bound in Corollary~\ref{cor3-5} and the upper bound in Proposition~\ref{pro3-10}, we have
	\begin{equation*}
		D\subset\{(1-2\delta)R\le r\le(1+2\delta)R\} \quad\text{for all sufficiently large }k.
	\end{equation*}
	Since $\delta$ can be arbitrarily small, the first limit in \eqref{3-44} follows.
\end{proof}

\section{Size and shape of the cross-section}\label{sec4}

By Proposition~\ref{pro3-11}, the radial width and the height of $D$ are both $o(R)$. We now bound both by a constant times $a$, the radius of a disk of area $|D|$. Recall the definitions of $a$ and the radial coordinate of the center of $D$,
\begin{align}
	a&=\left(\frac{|D|}{\pi}\right)^{1/2},\label{4-1}\\
	R_{\mathrm{center}}&=\frac1{|D|}\int_D r\,dr\,dz.\label{4-2}
\end{align}
All limits in this section are as $k\to\infty$, uniformly over centered finite-energy solutions. By \eqref{3-44}, the area of $D$ is $o(R^2)$ and the average of $r/R$ on $D$ tends to one. Hence $a=o(R)$ and $R_{\mathrm{center}}/R\to1$. Since $r/R\to1$ uniformly on $D$,
\begin{equation}\label{4-3}
	\mu(D)=15\pi R^{3/2}a^2(1+o(1)).
\end{equation}

\subsection{A bound on the height}\label{sec4-1}

By \eqref{3-44}, $r'/r$ stays in a fixed compact subinterval of $(0,\infty)$ on $D\times D$. Thus we can bound the vertical derivative of the potential by a constant times the integral of $1/|x-y|$ on $D$. By rearrangement, this integral is at most $2\pi a$.

\begin{lemma}\label{lem4-1}
	For all sufficiently large $k$, there is a universal $C$ such that
	\begin{equation}\label{4-4}
		\|F_z\|_{L^\infty(D)}\le Ca,\qquad
		\|F_+\|_{L^\infty(\R^5)}\le Ca^2,
	\end{equation}
	and
	\begin{equation}\label{4-11}
		\sup_D|z|\le Ca.
	\end{equation}
	If $\{z:(r,z)\in D\}=(-H(r),H(r))$, then
	\begin{equation}\label{4-6}
		\int_\R F_+(r,z)\,dz\le CaH(r)^2.
	\end{equation}
\end{lemma}

\begin{proof}
	Differentiating \eqref{2-3} in $z$, we have
	\begin{equation*}
		F_z(r,z)=v_z(r,z)
		=30r^{-5/2}\int_D
		(z-z')G'(\rho((r,z),(r',z')))(r')^{1/2}\,dr'\,dz'.
	\end{equation*}
	The kernel is locally integrable since $-G'(\rho)\le C/\rho$. For large $k$, \eqref{3-44} implies $R/2\le r,r'\le 2R$ on $D$, so
	\begin{equation}\label{4-7}
		\begin{split}
			|F_z(r,z)| &\le C\int_D\left(\frac{r'}{r}\right)^{3/2} \frac{|z-z'|}{(r-r')^2+(z-z')^2}\,dr'\,dz'\\
			&\le C\int_D\frac{dr'\,dz'} {\sqrt{(r-r')^2+(z-z')^2}} \le Ca.
		\end{split}
	\end{equation}
	Since $|D|=\pi a^2$, the rearrangement bound is
	\begin{equation*}
		\int_D\frac{dy}{|x-y|}
		\le\int_{B_a(0)}\frac{dy}{|y|}
		=2\pi\int_0^a\,dt=2\pi a.
	\end{equation*}

	Using \eqref{4-7} in the vertical integral in \eqref{3-16}, we have
	\begin{equation*}
		0<F(r,z)\le Ca\bigl(H(r)-|z|\bigr)
		\qquad(|z|<H(r)).
	\end{equation*}
	By Lemma~\ref{lem2-5},
	\begin{equation*}
		\|F_+\|_{L^\infty(\R^5)}
		\le Ca\sup_D|z|
		\le Ca\sqrt{\|F_+\|_{L^\infty(\R^5)}}.
	\end{equation*}
	Dividing by $\sqrt{\|F_+\|_{L^\infty(\R^5)}}$, we have $\|F_+\|_{L^\infty(\R^5)}\le Ca^2$. Lemma~\ref{lem2-5} then bounds the height by $Ca$. Integrating the bound for $F(r,z)$ over the section, we obtain
	\begin{equation*}
		\int_\R F_+(r,z)\,dz
		\le2Ca\int_0^{H(r)}\bigl(H(r)-z\bigr)\,dz
		=CaH(r)^2.
	\end{equation*}
\end{proof}

\subsection{The radial width of the core}\label{sec4-2}

Let $r_-$ and $r_+$ be the infimum and supremum of $r$ on $D$, and let $\ell$ be the larger of their distances from $R_{\mathrm{center}}$.
\begin{equation}\label{4-17}
	r_-=\inf_D r,\quad r_+=\sup_D r,\quad
	\ell=\max\{R_{\mathrm{center}}-r_-,r_+-R_{\mathrm{center}}\}=o(R).
\end{equation}
In the proof of Proposition~\ref{pro4-3}, we use \eqref{3-12} with $T$ near $r_-$ if $R_{\mathrm{center}}-r_-=\ell$ and near $r_+$ if $r_+-R_{\mathrm{center}}=\ell$. The following lemma bounds $R_{\mathrm{center}}-R$ so that we can estimate $B'$. It also shows that at least one third of $D$ lies on the other side of $r=T$, at distance at least $\ell/4$ from the line for the values of $T$ used below.

\begin{lemma}\label{lem4-2}
	Let $R_{\mathrm{center}}$ and $\ell$ be as in \eqref{4-2} and \eqref{4-17}. For all sufficiently large $k$,
	\begin{equation}\label{4-9}
		0\le R_{\mathrm{center}}-R\le C(Ra^2+\ell^2/R).
	\end{equation}
	Moreover,
	\begin{equation}\label{4-10}
		|D\cap\{r\le R_{\mathrm{center}}+\ell/2\}|\ge |D|/3,\qquad
		|D\cap\{r\ge R_{\mathrm{center}}-\ell/2\}|\ge |D|/3.
	\end{equation}
\end{lemma}

\begin{proof}
	By \eqref{3-9}, the expression $(r^{-3/2}\partial_r+(r')^{-3/2}\partial_{r'})G(\rho)$ is nonnegative and bounded by $C(r^{-5/2}+(r')^{-5/2})$. The factor $r^{3/2}(r')^{3/2}$ in \eqref{3-14} is bounded by $CR^3$ since $R/2\le r,r'\le 2R$ by \eqref{3-44}. Integrating over $D\times D$ and dividing by $|D|$, we obtain
	\begin{equation}\label{4-18}
		0\le \frac{1}{|D|}\int_D B'\,dr\,dz
		\le CR^{3-5/2}|D|=C\pi\sqrt R\,a^2 .
	\end{equation}
	For $R/2\le r\le 2R$, an interval containing $r_-,r_+$ and $R$ for large $k$,
	\begin{equation}\label{4-12}
		\begin{aligned}
			B'(R)&=0,\\
			B''(r)&=\frac34r^{-1/2}+\frac94R^2r^{-5/2},
			& cR^{-1/2}&\le B''(r)\le CR^{-1/2},\\
			B'''(r)&=-\frac38r^{-3/2}-\frac{45}{8}R^2r^{-7/2}<0,
			& |B'''(r)|&\le CR^{-3/2}.
		\end{aligned}
	\end{equation}
	By concavity of $B'$,
	\begin{equation*}
		0\le |D|^{-1}\int_D B'\,dr\,dz\le B'(R_{\mathrm{center}}).
	\end{equation*}
	Thus $R_{\mathrm{center}}\ge R$. By the bound on $B'''$, Taylor's formula at $R_{\mathrm{center}}$ implies
	\begin{equation*}
		B'(r)\ge B'(R_{\mathrm{center}})
		+B''(R_{\mathrm{center}})(r-R_{\mathrm{center}})
		-CR^{-3/2}(r-R_{\mathrm{center}})^2 .
	\end{equation*}
	The term linear in $r-R_{\mathrm{center}}$ integrates to zero by the definition of $R_{\mathrm{center}}$. Averaging over $D$ and using \eqref{4-18}, we have
	\begin{align*}
		B'(R_{\mathrm{center}}) &\le\frac1{|D|}\int_D B'\,dr\,dz +\frac C{R^{3/2}|D|}\int_D(r-R_{\mathrm{center}})^2\,dr\,dz\\
		&\le C\sqrt R\,a^2+CR^{-3/2}\ell^2 .
	\end{align*}
	On $[R,R_{\mathrm{center}}]$, the bound $B''(r)\ge cR^{-1/2}$ from \eqref{4-12} implies
	\begin{equation*}
		B'(R_{\mathrm{center}})=\int_R^{R_{\mathrm{center}}}B''(r)\,dr
		\ge cR^{-1/2}(R_{\mathrm{center}}-R).
	\end{equation*}
	This proves \eqref{4-9}.

	By the definition of $\ell$, $|r-R_{\mathrm{center}}|\le\ell$ on $D$. On $\{r>R_{\mathrm{center}}+\ell/2\}$ we have $r-R_{\mathrm{center}}\ge\ell/2$, and on the complement in $D$ we have $r-R_{\mathrm{center}}\ge-\ell$. Applying these bounds and the corresponding ones for $R_{\mathrm{center}}-r$ and integrating over $D$, we have
	\begin{align*}
		0=\int_D(r-R_{\mathrm{center}})\,dr\,dz &\ge \frac{3\ell}{2}|D\cap\{r>R_{\mathrm{center}}+\ell/2\}|-\ell|D|,\\
		0=\int_D(R_{\mathrm{center}}-r)\,dr\,dz &\ge \frac{3\ell}{2}|D\cap\{r<R_{\mathrm{center}}-\ell/2\}|-\ell|D|.
	\end{align*}
	Each of the two sets in the right-hand sides has area at most $2|D|/3$. Taking the complement in $D$ proves \eqref{4-10}.
\end{proof}

\begin{proposition}\label{pro4-3}
	For all sufficiently large $k$, the cross-section satisfies
	\begin{equation}\label{4-16}
		D\subset B_{Ca}((R_{\mathrm{center}},0)).
	\end{equation}
\end{proposition}

\begin{proof}
	By \eqref{4-11} it is enough to prove $\ell\le Ca$. We consider $\ell/a$ large enough that the height bound \eqref{4-11} implies
	\begin{equation}\label{4-14}
		|z-z'|\le\ell/8\qquad(x,y\in D).
	\end{equation}
	We first treat the case $\ell=r_+-R_{\mathrm{center}}$. Take $R_{\mathrm{center}}+3\ell/4<T<r_+$ and $D_0=D\cap\{r\le R_{\mathrm{center}}+\ell/2\}$. For $x=(r,z)\in D^+(T)$ and $y=(r',z')\in D_0$, we have $\ell/4\le r-r'\le2\ell$. By \eqref{4-14}, $(z-z')^2\le\ell^2/64\le(\ell/16)(r-r')$, and thus
	\begin{align*}
		r^2-(r')^2-(z-z')^2 &\ge(r+r'-\ell/16)(r-r')\ge cR(r-r'),\\
		\rho(x,y)&\le\frac{4\ell^2+\ell^2/64}{rr'} \le C(\ell/R)^2=o(1).
	\end{align*}
	By \eqref{3-10} and \eqref{2-4},
	\begin{equation}\label{4-15}
		-\partial_r G(\rho(x,y))
		=[-G'(\rho)]\frac{r^2-(r')^2-(z-z')^2}{r^2r'}
		\ge c\frac{r-r'}{(r-r')^2+(z-z')^2}\ge\frac c\ell .
	\end{equation}
	By \eqref{4-10}, $|D_0|\ge |D|/3$. Since $r\ge R/2$ on $D$, $\mu(D_0)\ge15(R/2)^{3/2}|D_0|\ge cR^{3/2}a^2$. Hence
	\begin{equation*}
		-\int_{D^+(T)}\int_{D_0}\partial_r G(\rho(x,y))\,d\mu(y)\,dx
		\ge\frac{cR^{3/2}a^2}{\ell}S^+(T).
	\end{equation*}
	Here $dx=dr\,dz$. We use \eqref{3-25} when $y\in D^+(T)$ and \eqref{3-26} when $R_{\mathrm{center}}+\ell/2<r'<T$. Since $15r^{3/2}\,dx=d\mu(x)$ on $D$, we have
	\begin{align*}
		\bigg[\int_{D^+(T)}\int_{D\setminus D_0}\partial_rG(\rho(x,y)) \,d\mu(y)\,dx\bigg]_+
		&\le CR^{-3/2}\mu(D\setminus D_0)\int_{D^+(T)}\frac{d\mu(x)}r\\
		&\le C\sqrt R\,a^2S^+(T).
	\end{align*}
	Since $r>T>R_{\mathrm{center}}\ge R$ on $D^+(T)$, $B'\ge0$ there and
	\begin{equation*}
		\int_{D^+(T)}(B'-V_r)\,dx
		\ge\left(\frac{cR^{3/2}a^2}{\ell}-C\sqrt R\,a^2\right)S^+(T).
	\end{equation*}
	If $\ell=R_{\mathrm{center}}-r_-$, take $r_-<T<R_{\mathrm{center}}-3\ell/4$ and $D_0=D\cap\{r\ge R_{\mathrm{center}}-\ell/2\}$. For $x\in D^-(T)$ and $y\in D_0$, $\ell/4\le r'-r\le2\ell$, so
	\begin{equation*}
		\partial_r G(\rho(x,y))
		=[-G'(\rho)]\frac{(r')^2-r^2+(z-z')^2}{r^2r'}
		\ge c\frac{r'-r}{(r-r')^2+(z-z')^2}\ge\frac c\ell.
	\end{equation*}
	The integral over $D^-(T)\times D^-(T)$ is nonnegative by \eqref{3-25} and $15\,dx=r^{-3/2}\,d\mu(x)$. For $y\in D^+(T)\setminus D_0$, $r'>r$ and thus $\partial_rG(\rho)>0$ by \eqref{3-10}. By \eqref{4-10}, $|D_0|\ge |D|/3$, just as in the first case. Using the integral over $D^-(T)\times D_0$, we obtain
	\begin{equation*}
		\int_{D^-(T)}V_r\,dx\ge\frac{c\mu(D_0)}\ell S^-(T)
		\ge\frac{cR^{3/2}a^2}\ell S^-(T).
	\end{equation*}
	By \eqref{4-9}, for $r<T$,
	\begin{equation*}
		r-R\le-\frac{3}{4}\ell+C(Ra^2+\ell^2/R)
		\le-\frac{1}{2}\ell+CRa^2 .
	\end{equation*}
	It follows that $[B'(r)]_+\le CR^{-1/2}[r-R]_+\le C\sqrt R\,a^2$. Hence
	\begin{equation*}
		\int_{D^-(T)}(V_r-B')\,dx
		\ge\left(\frac{cR^{3/2}a^2}{\ell}-C\sqrt R\,a^2\right)S^-(T).
	\end{equation*}
	Since $\ell/R\to0$, for large $k$ the term $C\sqrt R\,a^2$ is at most half of $cR^{3/2}a^2/\ell$. Write $S=S^+$ in the first case and $S=S^-$ in the second. Using \eqref{3-12}, we get the lower bound below on the chosen interval. The upper bound follows from \eqref{4-6} and $R/2\le T\le2R$. Thus
	\begin{equation*}
		CR^{3/2}aH(T)^2
		\ge T^{3/2}\int_\R F_+(T,z)\,dz
		\ge cR^{3/2}\frac{a^2}{\ell}S(T).
	\end{equation*}
	Dividing by $R^{3/2}a$, we have $H(T)^2\ge c(a/\ell)S(T)$. On this interval, $S>0$ and $|S'|=2H$ almost everywhere. By Lemma~\ref{lem3-2}, the area between its endpoints satisfies
	\begin{equation*}
		|D|\ge |S(T_1)-S(T_0)|
		\ge c\frac a\ell\left(\frac\ell4\right)^2
		=\frac{c}{16}a\ell,
	\end{equation*}
	where $T_0,T_1$ are the endpoints of the interval. Thus $\ell\le Ca$, and \eqref{4-11} bounds $|z|$ by $Ca$.
\end{proof}

For $x\in\partial D$ and $y\in D$, Proposition~\ref{pro4-3} implies $|x-y|\le Ca$. We can therefore use the expansion of $G$ in \eqref{2-4} in the formula for $v$ in \eqref{2-3}. Evaluating at the boundary, where $v=1+k/r^2$, determines the limit of $a^2\ln R$.

\begin{proposition}\label{pro4-4}
	As $k\to\infty$,
	\begin{equation}\label{4-46}
		a^2\ln R\to\frac8{15}.
	\end{equation}
	In particular, $a\to0$ and $aR_{\mathrm{center}}\to\infty$.
\end{proposition}

\begin{proof}
	By Proposition~\ref{pro3-11}, $a=o(R)$ and $R_{\mathrm{center}}/R\to1$. Take $x_0=(r_0,z_0)\in\partial D$ and write $y=(r',z')$ in \eqref{2-3}. Then
	\begin{equation*}
		4+o(1)=1+\frac{k}{r_0^2}=v(x_0)
		=15\int_D\left(\frac{r'}{r_0}\right)^{3/2}
		G(\rho(x_0,y))\,dr'\,dz' .
	\end{equation*}
	By Proposition~\ref{pro4-3}, $|x_0-y|\le Ca$. By \eqref{3-44}, $r_0/R\to1$ and $r'/R\to1$ uniformly for $x_0\in\partial D$ and $y\in D$. Consequently,
	\begin{align*}
		\int_D\left(\frac{r'}{r_0}\right)^{3/2}\,dr'\,dz'&=(1+o(1))\pi a^2,\quad \ln\frac{r_0r'}{R^2}=o(1),\\
		0<\rho(x_0,y)&\le C(a/R)^2,\\
		\sup_{y\in D}\rho(x_0,y)|\ln\rho(x_0,y)|&\le C(a/R)^2\bigl(1+\ln(R/a)\bigr)\to0.
	\end{align*}
	Expanding $G$ by \eqref{2-4},
	\begin{align*}
		G(\rho(x_0,y)) ={}&\frac1{2\pi}\ln\frac Ra +\frac1{2\pi}\ln\frac a{|x_0-y|} +\frac1{4\pi}\ln\frac{r_0r'}{R^2}+c_0\\
		&\quad+O\bigl(\rho(x_0,y)|\ln\rho(x_0,y)|\bigr).
	\end{align*}
	By $(r'/r_0)^{3/2}\le C$ and rearrangement,
	\begin{align*}
		\int_D\left(\frac{r'}{r_0}\right)^{3/2} \left(\ln\frac a{|x_0-y|}\right)_+\,dr'\,dz' &\le C\int_{B_a(0)}\ln\frac a{|y|}\,dy\\
		&=2\pi C\int_0^a t\ln\frac at\,dt =\frac{\pi C}{2}a^2.
	\end{align*}
	Since $|x_0-y|\le Ca$, the negative part of $\ln(a/|x_0-y|)$ is uniformly bounded on $D$. Its integral is therefore $O(a^2)$. The terms $\ln(r_0r'/R^2)$, $c_0$ and $O(\rho|\ln\rho|)$ are also bounded uniformly on $D$. Substituting into the formula for $v(x_0)$, we obtain
	\begin{equation}\label{4-47}
		\begin{split}
			4+o(1) &=\left(\frac{15}{2}+o(1)\right)a^2\ln\frac Ra+O(a^2)\\
			&=a^2\ln\frac Ra \left(\frac{15}{2}+o(1) +O\left(\frac1{\ln(R/a)}\right)\right).
		\end{split}
	\end{equation}
	Since $\ln(R/a)\to\infty$, $a^2\ln(R/a)\to8/15$ and $a\to0$. Then
	\begin{equation*}
		a^2\ln R
		=a^2\ln(R/a)-a^2\ln(1/a)\to\frac8{15}.
	\end{equation*}
	Since $(aR)^2=(a^2\ln R)R^2/\ln R\to\infty$ and $R_{\mathrm{center}}/R\to1$, we have $aR_{\mathrm{center}}\to\infty$. For the limit involving $R_{\mathrm{center}}$ in Theorem~\ref{thm1-2}, we have
	\begin{equation*}
		a^2\ln(R_{\mathrm{center}}/a)
		=a^2\ln(R/a)+a^2\ln(R_{\mathrm{center}}/R)
		\to\frac8{15}.
	\end{equation*}
\end{proof}

\subsection{The disk limit}\label{sec4-3}

For the comparison in Section~\ref{sec5}, we choose a length $b$ depending only on $k$. By Proposition~\ref{pro4-4}, $b/a\to1$ for every centered solution. After translating the center of $D$ to the origin and scaling by $b$, we write
\begin{equation}\label{4-29}
	\begin{aligned}
		b&:=\left(\frac{8}{15\ln R}\right)^{1/2},\\
		X&:=\left(\frac{r-R_{\mathrm{center}}}{b},\frac{z}{b}\right),&
		\Omega&:=\frac{D-(R_{\mathrm{center}},0)}{b},\\
		r(X)&:=R_{\mathrm{center}}+bX_1,&
		\tilde F(X)&:=\frac{F(r(X),bX_2)}{15b^2}.
	\end{aligned}
\end{equation}
By Proposition~\ref{pro4-4},
\begin{equation}\label{4-19}
	\frac{a^2}{b^2}=\frac{15}{8}a^2\ln R\to1,\qquad
	\frac{R_{\mathrm{center}}}{R}\to1.
\end{equation}
By Proposition~\ref{pro4-3}, $\Omega$ lies in a fixed ball $B_M(0)$ in $\R^2$, with $M>2$. The area and first moment of $\Omega$ satisfy
\begin{equation}\label{4-20}
	|\Omega|=\pi(a/b)^2\to\pi,\qquad
	\int_\Omega X\,dX=0.
\end{equation}
The first component of the moment is zero by the definition of $R_{\mathrm{center}}$, and the second by evenness in $z$. The half-plane $r(X)>0$ contains each fixed ball about the origin for large $k$. In this half-plane $\Omega=\{\tilde F>0\}$.

We first expand the kernel in the rescaled variables. The leading terms are the logarithmic kernel in $\R^2$ and a constant. To compare two solutions with different centers in Section~\ref{sec5}, we also estimate the derivative of the remainder with respect to $R_{\mathrm{center}}$ on $[R/2,2R]$.

\begin{lemma}\label{lem4-5}
	Fix $L>0$ and take $R$ sufficiently large. For each $R_{\mathrm{center}}\in[R/2,2R]$, use $b$ and $r(X)$ from \eqref{4-29}. For $X,Y\in B_L(0)$ with $X\ne Y$, evaluate $\rho$ at $x=(r(X),bX_2)$ and $y=(r(Y),bY_2)$. Then
	\begin{equation}\label{4-21}
		\left(\frac{r(Y)}{r(X)}\right)^{3/2}G(\rho)
		=\frac1{2\pi}\ln\frac{R_{\mathrm{center}}}{b}+c_0
		+\frac1{2\pi}\ln\frac1{|X-Y|}+K(X,Y),
	\end{equation}
	where $c_0$ is the constant in \eqref{2-4}, and
	\begin{equation}\label{4-22}
		|K(X,Y)|+R|\partial_{R_{\mathrm{center}}}K(X,Y)|
		\le C_L\frac{b\ln R}{R}.
	\end{equation}
	In taking this derivative we keep $R,b,X,Y$ fixed.
\end{lemma}

\begin{proof}
	By the definition of $\rho$ in \eqref{2-9},
	\begin{equation}\label{4-23}
		\rho=\frac{b^2|X-Y|^2}{r(X)r(Y)},\qquad
		\partial_{R_{\mathrm{center}}}\rho
		=-\rho\left(\frac1{r(X)}+\frac1{r(Y)}\right).
	\end{equation}
	The logarithmic term in \eqref{2-4} is
	\begin{equation*}
		\frac1{4\pi}\ln\frac1\rho
		=\frac1{2\pi}\ln\frac{R_{\mathrm{center}}}{b}
		+\frac1{2\pi}\ln\frac1{|X-Y|}
		+\frac1{4\pi}\ln\frac{r(X)r(Y)}{R_{\mathrm{center}}^2}.
	\end{equation*}
	Substituting into \eqref{4-21}, we have
	\begin{align*}
		K(X,Y)&=\left[\left(\frac{r(Y)}{r(X)}\right)^{3/2}-1\right]
		\left(\frac1{2\pi}\ln\frac{R_{\mathrm{center}}}{b}+c_0+\frac1{2\pi}\ln\frac1{|X-Y|}\right)\\
		&\quad+\left(\frac{r(Y)}{r(X)}\right)^{3/2}
		\left[\frac1{4\pi}\ln\frac{r(X)r(Y)}{R_{\mathrm{center}}^2}
		+G(\rho)-\frac1{4\pi}\ln\frac1\rho-c_0\right].
	\end{align*}
	For the ratio $r(Y)/r(X)$, we have
	\begin{equation*}
		\frac{r(Y)}{r(X)}-1=\frac{b(Y_1-X_1)}{r(X)},\qquad
		\partial_{R_{\mathrm{center}}}\!\left(\frac{r(Y)}{r(X)}\right)
		=-\frac{b(Y_1-X_1)}{r(X)^2}.
	\end{equation*}
	Since $r(X),r(Y)\in[R/4,3R]$ for large $R$, the mean value theorem and the chain rule imply
	\begin{equation*}
		\left|\left(\frac{r(Y)}{r(X)}\right)^{3/2}-1\right|
		+R\left|\partial_{R_{\mathrm{center}}}\left(\frac{r(Y)}{r(X)}\right)^{3/2}\right|
		\le C_L\frac bR|X-Y|.
	\end{equation*}
	We have $\ln(R_{\mathrm{center}}/b)=O(\ln R)$ and $\partial_{R_{\mathrm{center}}}\ln(R_{\mathrm{center}}/b)=1/R_{\mathrm{center}}$. The first product in $K$ and $R$ times its derivative are therefore bounded by
	\begin{equation*}
		C_L\frac bR |X-Y|\bigl(1+\ln R+\left|\ln\frac1{|X-Y|}\right|\bigr)
		\le C_L\frac{b\ln R}{R}.
	\end{equation*}
	For the second product in $K$, we use
	\begin{align*}
		\ln\frac{r(X)}{R_{\mathrm{center}}}
		&=\ln\left(1+\frac{bX_1}{R_{\mathrm{center}}}\right)=O(b/R),\\
		\partial_{R_{\mathrm{center}}}\ln\frac{r(X)}{R_{\mathrm{center}}}
		&=-\frac{bX_1}{R_{\mathrm{center}}\,r(X)}=O(b/R^2),
	\end{align*}
	and the same bounds at $Y$. By \eqref{2-4},
	\begin{equation*}
		\left|G(\rho)-\frac1{4\pi}\ln\frac1\rho-c_0\right|
		+\rho\left|G'(\rho)+\frac1{4\pi\rho}\right|
		\le C\rho(1+|\ln\rho|).
	\end{equation*}
	By \eqref{4-23}, $|\partial_{R_{\mathrm{center}}}\rho|\le C\rho/R$ and $\rho\le C_L(b/R)^2$. In particular,
	\begin{equation*}
		R\left|\partial_{R_{\mathrm{center}}}\left(G(\rho)-\frac1{4\pi}\ln\frac1\rho-c_0\right)\right|
		=R\left|G'(\rho)+\frac{1}{4\pi\rho}\right||\partial_{R_{\mathrm{center}}}\rho|
		\le C\rho(1+|\ln\rho|).
	\end{equation*}
	The bound $\rho\le C_L(b/R)^2$ and monotonicity of $t(1+|\ln t|)$ for $0<t<1$ imply $\rho(1+|\ln\rho|)\le C_L(b/R)^2(1+\ln(R/b))$. Thus the second product in $K$ and $R$ times its derivative are bounded by $C_L(b/R+(b/R)^2(1+\ln(R/b)))$. This is bounded by $C_L b\ln R/R$, as is the contribution of the first product.
\end{proof}

We integrate \eqref{4-21} over $\Omega$ to expand $\tilde F$. For two solutions with the same $k$, we use \eqref{4-22} and the mean value theorem with respect to $R_{\mathrm{center}}$.

\begin{corollary}\label{lem4-6}
	For any fixed $L>M$ and all sufficiently large $k$, the potential in \eqref{4-29} satisfies, on $B_L(0)$,
	\begin{equation}\label{4-24}
		\begin{split}
			\tilde F(X)&=c(k,v)+\frac1{2\pi}\int_\Omega\ln\frac1{|X-Y|}\,dY+E(X),\\
			c(k,v)&=\left(\frac1{2\pi}\ln\frac{R_{\mathrm{center}}}{b}+c_0\right)|\Omega|
			-\frac{1+k/R_{\mathrm{center}}^2}{15b^2}.
		\end{split}
	\end{equation}
	The error satisfies
	\begin{equation}\label{4-25}
		\sup_{X\in B_L(0)}|E(X)|\le C_L\frac{b\ln R}{R}.
	\end{equation}
	If two solutions have the same $k$, their corresponding errors satisfy
	\begin{equation}\label{4-26}
		\sup_{X\in B_L(0)}|E_1(X)-E_2(X)|
		\le C_L\frac{b\ln R}{R}
		\left(|\Omega_1\bigtriangleup\Omega_2|
		+\frac{|R_{\mathrm{center},1}-R_{\mathrm{center},2}|}{R}\right).
	\end{equation}
\end{corollary}

\begin{proof}
	In \eqref{2-3}, the change of variables $dr'\,dz'=b^2\,dY$ cancels the factor $b^{-2}$ in $\tilde F$. Hence
	\begin{equation*}
		\tilde F(X)=\int_\Omega\left(\frac{r(Y)}{r(X)}\right)^{3/2}G(\rho)\,dY
		-\frac{1+k/r(X)^2}{15b^2}.
	\end{equation*}
	By \eqref{4-19}, $R_{\mathrm{center}}\in[R/2,2R]$ for large $k$. Since $\Omega\subset B_M(0)$, we apply \eqref{4-21} inside the integral to obtain \eqref{4-24} with
	\begin{equation*}
		E(X)=\int_\Omega K(X,Y)\,dY
		-\frac{k}{15b^2}\bigl(r(X)^{-2}-R_{\mathrm{center}}^{-2}\bigr).
	\end{equation*}
	For $R_{\mathrm{center}}\in[R/2,2R]$, we have
	\begin{align*}
		r(X)^{-2}-R_{\mathrm{center}}^{-2}
		&=-bX_1\frac{r(X)+R_{\mathrm{center}}}{r(X)^2 R_{\mathrm{center}}^2},\\
		\partial_{R_{\mathrm{center}}}\bigl(r(X)^{-2}-R_{\mathrm{center}}^{-2}\bigr)
		&=-2\bigl(r(X)^{-3}-R_{\mathrm{center}}^{-3}\bigr).
	\end{align*}
	The absolute value of the first line is at most $C_Lb/R^3$. For the second line, the mean value theorem applied to $s^{-3}$ and $|r(X)-R_{\mathrm{center}}|=b|X_1|$ yield the bound $C_Lb/R^4$. After multiplication by $k/(15b^2)=R^2/(5b^2)$, the term in $E$ containing $k$ is bounded by $C_L/(bR)$ and its derivative by $C_L/(bR^2)$. Using \eqref{4-22}, $|\Omega|\le \pi M^2$ and $b^2\ln R=8/15$, we obtain
	\begin{equation*}
		\|E\|_{L^\infty(B_L(0))}
		\le C_L\left(|\Omega|\frac{b\ln R}{R}+\frac1{bR}\right)
		\le C_L\frac{b\ln R}{R}.
	\end{equation*}

	For two solutions at the same $k$, write $K_i,r_i$ for the functions corresponding to $R_{\mathrm{center},i}$. Subtracting the two errors, we have
	\begin{align*}
		E_1(X)-E_2(X)
		&=\int_{\Omega_1\setminus\Omega_2}K_1(X,Y)\,dY-\int_{\Omega_2\setminus\Omega_1}K_1(X,Y)\,dY\\
		&\quad+\int_{\Omega_2}\bigl(K_1(X,Y)-K_2(X,Y)\bigr)\,dY\\
		&\quad-\frac{k}{15b^2}\left[r_1(X)^{-2}-R_{\mathrm{center},1}^{-2}
		-r_2(X)^{-2}+R_{\mathrm{center},2}^{-2}\right].
	\end{align*}
	By \eqref{4-22}, we bound the first two integrals by $C_L(b\ln R/R)|\Omega_1\bigtriangleup\Omega_2|$. For $Y\in\Omega_2$, the derivative bound in \eqref{4-22} and the mean value theorem imply
	\begin{equation*}
		|K_1(X,Y)-K_2(X,Y)|
		\le C_L\frac{b\ln R}{R}\frac{|R_{\mathrm{center},1}-R_{\mathrm{center},2}|}{R}.
	\end{equation*}
	This bound is uniform for $Y\in\Omega_2$, so we multiply it by $|\Omega_2|$ to estimate the integral. The difference of the terms involving $k$ is at most $C_L|R_{\mathrm{center},1}-R_{\mathrm{center},2}|/(bR^2)$ by the derivative bound above. Since $1/(bR^2)=(15/8)b\ln R/R^2$ and $|\Omega_2|\le\pi M^2$, summing the estimates proves \eqref{4-26}.
\end{proof}

Friedman and Turkington \cite{FT} obtained the disk limit for their constrained maximizers. For centered finite-energy solutions, we obtain the limiting potential from \eqref{4-24} and \eqref{4-25}. We then use this limit to describe the boundary of $\Omega$.

\begin{proposition}\label{pro4-7}
	For every $0<\alpha<1$,
	\begin{equation*}
		\overline\Omega\to\overline{B_1(0)}
		\quad\text{in Hausdorff distance},\qquad
		\tilde F\to\tilde F_*
		\quad\text{in }C^{1,\alpha}_{\mathrm{loc}}(\R^2),
	\end{equation*}
	where
	\begin{equation}\label{4-27}
		\tilde F_*(X)=
		\begin{cases}
			(1-|X|^2)/4,& |X|\le1,\\
			\frac{1}{2}\ln\frac1{|X|},& |X|\ge1 .
		\end{cases}
	\end{equation}
	For all sufficiently large $k$, $\Omega$ has the form
	\begin{equation}\label{4-28}
		\Omega=\{s e(\theta):0\le s<s_0(\theta)\},\qquad
		e(\theta)=(\cos\theta,\sin\theta),\qquad
		s_0\to1\quad\text{in }C^{1,\alpha}(\mathbb T).
	\end{equation}
	Consequently, $\Omega$ is simply connected. The domain $(D-(R_{\mathrm{center}},0))/a$ also converges to the unit disk in Hausdorff distance, and its radius in polar coordinates converges to one in $C^{1,\alpha}(\mathbb T)$.
\end{proposition}

\begin{proof}
	For $L>M$, by translation of the integration variable and polar coordinates,
	\begin{equation*}
		\sup_{X\in B_{2L}(0)}\int_{B_M(0)}\left|\ln\frac1{|X-Y|}\right|\,dY
		\le 2\pi\int_0^{M+2L}t|\ln t|\,dt<\infty.
	\end{equation*}
	Take $X_0\in\partial\Omega$. Since $\tilde F(X_0)=0$, \eqref{4-24} at $X_0$ and \eqref{4-25} imply
	\begin{equation*}
		|c(k,v)|\le\frac1{2\pi}\int_{B_M(0)}\left|\ln\frac1{|X_0-Y|}\right|\,dY
		+|E(X_0)|\le C_M.
	\end{equation*}
	Using the same bound on the integral of the logarithm for $X\in B_{2L}(0)$, we have
	\begin{equation*}
		\|\tilde F\|_{L^\infty(B_{2L}(0))}
		\le |c(k,v)|+\frac1{2\pi}\sup_{X\in B_{2L}(0)}\int_{B_M(0)}
		\left|\ln\frac1{|X-Y|}\right|\,dY+\|E\|_{L^\infty(B_{2L}(0))}
		\le C_L.
	\end{equation*}
	Since $(\partial_r^2+(3/r)\partial_r)r^{-2}=0$ for $r>0$, \eqref{1-1} and \eqref{4-29} imply
	\begin{equation*}
		\partial_{X_1}\tilde F=\frac{F_r}{15b},\qquad
		\Delta_X\tilde F=\frac{F_{rr}+F_{zz}}{15},\qquad
		-F_{rr}-F_{zz}-\frac3r F_r=15\boldsymbol 1_D .
	\end{equation*}
	Substituting,
	\begin{equation}\label{4-31}
		-\Delta_X\tilde F-\frac{3b}{r(X)}\partial_{X_1}\tilde F
		=\boldsymbol 1_\Omega .
	\end{equation}
	Since $3b/r(X)\le C_Lb/R$ on $B_{2L}(0)$, interior estimates give, for $1<p<\infty$,
	\begin{equation*}
		\|\tilde F\|_{W^{2,p}(B_L(0))}
		\le C_{L,p}\bigl(
		\|\tilde F\|_{L^p(B_{2L}(0))}+
		\|\boldsymbol 1_\Omega\|_{L^p(B_{2L}(0))}\bigr)\le C_{L,p}.
	\end{equation*}
	Given solutions $(k_n,v_n)$ with $k_n\to\infty$, take a subsequence such that
	\begin{align*}
		\tilde F_n&\to\tilde F_\infty \quad\text{in }C^{1,\alpha}_{\mathrm{loc}}(\R^2)\quad(0<\alpha<1),\\
		c(k_n,v_n)&\to c_\infty,\quad \boldsymbol 1_{\Omega_n}\rightharpoonup^* g\quad\text{in }L^\infty(\R^2).
	\end{align*}
	The interior bounds also imply $\tilde F_\infty\in W^{2,p}_{\mathrm{loc}}(\R^2)$ for every finite $p>1$. Since $\Omega_n\subset B_M(0)$, $g=0$ outside $\overline{B_M(0)}$. For each fixed $X$, the function $\boldsymbol 1_{B_M(0)}(Y)\ln(1/|X-Y|)$ is integrable in $Y$. We use this function to test the weak-* convergence of $\boldsymbol 1_{\Omega_n}$. Passing to the limit in \eqref{4-24} and using $E_n(X)\to0$, we have
	\begin{equation}\label{4-30}
		\tilde F_\infty(X)=c_\infty+\frac1{2\pi}\int_{B_M(0)}\ln\frac1{|X-Y|}g(Y)\,dY .
	\end{equation}
	The logarithmic kernel in \eqref{4-30} is the fundamental solution of $-\Delta_X$, so $-\Delta_X\tilde F_\infty=g$.

	On each compact subset of $\{\tilde F_\infty>0\}$, uniform convergence implies $\tilde F_n>0$ for large $n$. By weak-* convergence of $\boldsymbol 1_{\Omega_n}$, we have $g=1$ almost everywhere on $\{\tilde F_\infty>0\}$. The same argument shows that $g=0$ almost everywhere on $\{\tilde F_\infty<0\}$.

	It remains to show that $g=0$ almost everywhere on $\{\tilde F_\infty=0\}$. By the equation $g=-\Delta_X\tilde F_\infty$, we can argue as in \cite{Qin}. By Theorem~4.4(iv) in \cite{EG}, the weak gradient of a Sobolev function vanishes almost everywhere where the function is constant. Thus $\partial_j\tilde F_\infty=0$ almost everywhere on $\{\tilde F_\infty=0\}$ for $j=1,2$. Applying the same property to $\partial_j\tilde F_\infty\in W^{1,p}_{\mathrm{loc}}(\R^2)$, we deduce that $\nabla\partial_j\tilde F_\infty=0$ almost everywhere on $\{\partial_j\tilde F_\infty=0\}$. In particular,
	\begin{equation*}
		-\Delta_X\tilde F_\infty=0
		\qquad\text{almost everywhere on }\{\tilde F_\infty=0\}.
	\end{equation*}
	We conclude that $g=\boldsymbol 1_{\Omega_\infty}$ for $\Omega_\infty=\{\tilde F_\infty>0\}$. Passing to the limit in \eqref{4-20} with the tests $\boldsymbol 1_{B_M(0)}$ and $X_j\boldsymbol 1_{B_M(0)}$, we obtain
	\begin{equation*}
		|\Omega_\infty|=\pi,\qquad \int_{\Omega_\infty}X\,dX=0 .
	\end{equation*}
	The set $\Omega_\infty$ is bounded and open. By continuity, $\tilde F_\infty=0$ on $\partial\Omega_\infty$, so \eqref{4-30} shows that the logarithmic potential of $\Omega_\infty$ is constant on its boundary. By Fraenkel's characterization of bounded open sets with constant logarithmic potential on their boundary \cite{Fb} (see Theorem~1 in \cite{Rei}), $\Omega_\infty$ is a disk. Its area is $\pi$ and its first moment is zero, so $\Omega_\infty=B_1(0)$.

	With $\Omega_\infty=B_1(0)$ in \eqref{4-30}, the potential is radial and vanishes at $s=|X|=1$. Integrating $-\Delta_X\tilde F_\infty=\boldsymbol 1_{B_1(0)}$ over $B_s(0)$, we obtain
	\begin{equation*}
		-2\pi s\,\frac{d}{ds}\tilde F_\infty(s e(\theta))
		=\int_{B_s(0)}\boldsymbol 1_{B_1(0)}\,dX=\pi\min\{s^2,1\}.
	\end{equation*}
	The resulting radial derivative is $-s/2$ for $s<1$ and $-1/(2s)$ for $s>1$. Integrating from $s=1$, we obtain \eqref{4-27}. The sequence of solutions with $k_n\to\infty$ was arbitrary and its convergent subsequence has limit $\tilde F_*$. It follows that $\tilde F\to\tilde F_*$ in $C^{1,\alpha}_{\mathrm{loc}}(\R^2)$ as $k\to\infty$.

	To describe $\partial\Omega$ in polar coordinates, we argue as in \cite{Qin}. Fix $0<\eta<1/4$. By \eqref{4-27} and uniform convergence on $\overline{B_M(0)}$, for all sufficiently large $k$ we have
	\begin{equation*}
		\tilde F>0\quad\text{on }\overline{B_{1-\eta}(0)},\qquad
		\tilde F<0\quad\text{for }1+\eta\le|X|\le M.
	\end{equation*}
	The two formulas $-s/2$ and $-1/(2s)$ for $\partial_s\tilde F_*$ are bounded above by $-3/8$ on $1-\eta\le s\le1+\eta$. By $C^1$ convergence on $\overline{B_2(0)}$,
	\begin{equation*}
		\frac{d}{ds}\tilde F(se(\theta))=\nabla_X\tilde F(se(\theta))\cdot e(\theta)
		\le-\frac14
		\qquad(1-\eta\le s\le1+\eta).
	\end{equation*}
	For fixed $\theta$, the function $\tilde F(se(\theta))$ is positive at $s=1-\eta$, negative at $s=1+\eta$, and strictly decreasing between them. Thus it vanishes at a unique $s_0(\theta)\in(1-\eta,1+\eta)$. It is positive for $0\le s<s_0(\theta)$ and negative for $s_0(\theta)<s\le M$. Since $\Omega\subset B_M(0)$ and $\Omega=\{\tilde F>0\}$, we obtain \eqref{4-28}. The inclusions $B_{1-\eta}(0)\subset\Omega\subset B_{1+\eta}(0)$ for every $\eta$ imply Hausdorff convergence and $\|s_0-1\|_{L^\infty(\mathbb T)}\to0$.

	By the implicit function theorem applied to $\tilde F(s_0(\theta)e(\theta))=0$, we have $s_0\in C^{1,\alpha}(\mathbb T)$. Since $\tilde F_*$ is radial, $\nabla_X\tilde F_*(s_0(\theta)e(\theta))\cdot e'(\theta)=0$. Differentiating the equation for $s_0$ and using this identity, we get
	\begin{equation*}
		s_0'(\theta)
		=-\frac{s_0(\theta)\nabla_X(\tilde F-\tilde F_*)(s_0(\theta)e(\theta))\cdot e'(\theta)}
		{\nabla_X\tilde F(s_0(\theta)e(\theta))\cdot e(\theta)}.
	\end{equation*}
	The denominator is bounded away from zero by $1/4$. This formula bounds $\|s_0'\|_{L^\infty(\mathbb T)}$ uniformly. Thus $\theta\mapsto s_0(\theta)e(\theta)$ is uniformly Lipschitz, and
	\begin{equation*}
		\bigl\|\nabla_X(\tilde F-\tilde F_*)(s_0(\cdot)e(\cdot))\bigr\|_{C^\alpha(\mathbb T)}
		\le C\|\tilde F-\tilde F_*\|_{C^{1,\alpha}(B_2(0))}.
	\end{equation*}
	The same argument bounds the $C^\alpha(\mathbb T)$ norm of the denominator. Since its absolute value is at least $1/4$, the reciprocal is also bounded in $C^\alpha(\mathbb T)$. Hence $\|s_0'\|_{C^\alpha(\mathbb T)}\to0$. Thus $s_0\to1$ in $C^{1,\alpha}(\mathbb T)$. Every segment from the origin to a point of $\Omega$ lies in $\Omega$, so $\Omega$ is simply connected. After scaling by $a$ in place of $b$, the radius as a function of $\theta$ is $(b/a)s_0$. By \eqref{4-19}, it also converges to one in $C^{1,\alpha}(\mathbb T)$.
\end{proof}

\section{Uniqueness at large flux}\label{sec5}

Let $(k,v_i)$, $i=1,2$, be two centered solutions at the same flux. In \eqref{4-29}, they have the same $R$ and $b$, and $r_i(X)=R_{\mathrm{center},i}+bX_1$ for $i=1,2$. By Proposition~\ref{pro4-7}, the two rescaled cross-sections $\Omega_i$ have radii $s_i$ in polar coordinates. By symmetry with respect to $X_2=0$, each $s_i$ is even in $\theta$. Set
\begin{equation*}
	h:=s_1-s_2,\qquad
	d:=\frac{|R_{\mathrm{center},1}-R_{\mathrm{center},2}|}{R}+\|h\|_{L^2(\mathbb T)} .
\end{equation*}
Fix $0<\delta<1/4$ and take $k$ large. By \eqref{4-19} and Proposition~\ref{pro4-7},
\begin{equation*}
	\left|\frac{R_{\mathrm{center},i}}{R}-1\right|+\|s_i-1\|_{C^1(\mathbb T)}
	+\|\tilde F_i-\tilde F_*\|_{C^1(B_2(0))}<\delta,\qquad i=1,2.
\end{equation*}
The constants $C$ below are independent of $\delta$ and $k$. In polar coordinates,
\begin{equation}\label{5-1}
	\begin{split}
		|\Omega_1\bigtriangleup\Omega_2|&=\frac12\int_{\mathbb T}|s_1(\theta)^2-s_2(\theta)^2|\,d\theta
		\le C\|h\|_{L^2(\mathbb T)},\\
		|(\Omega_1\times\Omega_1)\bigtriangleup(\Omega_2\times\Omega_2)|
		&\le 2|\Omega_1\cup\Omega_2|\,|\Omega_1\bigtriangleup\Omega_2|.
	\end{split}
\end{equation}
For the second inequality, if $(X,Y)$ belongs to the symmetric difference of $\Omega_1\times\Omega_1$ and $\Omega_2\times\Omega_2$, at least one of $X$ and $Y$ belongs to $\Omega_1\bigtriangleup\Omega_2$. The other belongs to $\Omega_1\cup\Omega_2$. We also use $\bigl||\Omega_1|-|\Omega_2|\bigr|\le|\Omega_1\bigtriangleup\Omega_2|$.

We first estimate the difference of the centers by subtracting the radial identity \eqref{3-14} for the two solutions.

\begin{lemma}\label{lem5-1}
	For two centered solutions with the same sufficiently large $k$,
	\begin{equation}\label{5-2}
		\frac{|R_{\mathrm{center},1}-R_{\mathrm{center},2}|}{R}\le Cb^2d .
	\end{equation}
\end{lemma}

\begin{proof}
	Apply \eqref{3-14} to either solution and divide by $\sqrt R\,|D|$. Since $dx=b^2\,dX$, $dy=b^2\,dY$ and $|D|=b^2|\Omega|$, the double integral has the factor $b^2/|\Omega|$ after the change of variables. We obtain
	\begin{equation}\label{5-3}
		\begin{split}
			&\frac{1}{|\Omega|}\int_\Omega\frac{B'(r(X))}{\sqrt R}\,dX\\
			&\quad=\frac{15b^2}{2|\Omega|}\iint_{\Omega\times\Omega}
			\frac{r(X)^{3/2}r(Y)^{3/2}}{\sqrt R}
			\left(r(X)^{-3/2}\partial_r+r(Y)^{-3/2}\partial_{r'}\right)
			G(\rho(x,y))\,dX\,dY ,
		\end{split}
	\end{equation}
	Here we evaluate the derivatives $\partial_r$ and $\partial_{r'}$ from \eqref{3-9} at $x=(r(X),bX_2)$ and $y=(r(Y),bY_2)$. Substituting \eqref{3-9} and using $r(X)-r(Y)=b(X_1-Y_1)$ and $|x-y|=b|X-Y|$, we write the integrand on the right-hand side as
	\begin{align*}
		&[-\rho G'(\rho)]
		\left(\frac{r(X)}R\right)^{3/2}\left(\frac{r(Y)}R\right)^{3/2}
		\bigg[\left(\frac{r(X)}R\right)^{-5/2}+\left(\frac{r(Y)}R\right)^{-5/2}\\
		&\qquad+3\frac{(X_1-Y_1)^2}{|X-Y|^2}
		\int_0^1\left(\frac{(1-\tau)r(Y)+\tau r(X)}R\right)^{-5/2}\,d\tau\bigg].
	\end{align*}
	The quotient $(X_1-Y_1)^2/|X-Y|^2$ is at most one and is independent of the center. The ratios $r(X)/R$ and $r(Y)/R$ lie in a compact subinterval of $(0,\infty)$, and their derivatives with respect to $R_{\mathrm{center}}$ equal $1/R$. To differentiate $\rho G'(\rho)$, we use \eqref{2-2} to obtain
	\begin{equation*}
		G''(t)=\frac{15}{32\pi^2}\int_{\mathbb S^3}
		\bigl(t+2(1-\omega_1)\bigr)^{-7/2}\,d\omega
		\le\frac{5}{2t}[-G'(t)].
	\end{equation*}
	Indeed the denominator of the integrand for $G''$ has one extra factor $t+2(1-\omega_1)\ge t$ compared with the integrand for $-G'$. By \eqref{2-5} and \eqref{4-23}, we have
	\begin{equation*}
		\bigl|\partial_{R_{\mathrm{center}}}\bigl(\rho G'(\rho)\bigr)\bigr|
		\le \bigl(|G'(\rho)|+\rho|G''(\rho)|\bigr)|\partial_{R_{\mathrm{center}}}\rho|
		\le C\rho^{-1}|\partial_{R_{\mathrm{center}}}\rho|\le \frac CR.
	\end{equation*}
	Thus the integrand in \eqref{5-3} is uniformly bounded and its derivative with respect to $R_{\mathrm{center}}$ is bounded by $C/R$ for $X,Y\in B_M(0)$ and $R/2\le R_{\mathrm{center}}\le2R$. Let $f_i(X,Y)$ be this integrand evaluated at $R_{\mathrm{center},i}$. We compare the two integrals by writing
	\begin{align*}
		&\frac1{|\Omega_1|}\iint_{\Omega_1\times\Omega_1} f_1(X,Y)\,dX\,dY
		-\frac1{|\Omega_2|}\iint_{\Omega_2\times\Omega_2} f_2(X,Y)\,dX\,dY\\
		&\quad=\frac1{|\Omega_1|}
		\left(\iint_{\Omega_1\times\Omega_1} f_1\,dX\,dY
		-\iint_{\Omega_2\times\Omega_2} f_1\,dX\,dY\right)\\
		&\qquad+\frac1{|\Omega_1|}\iint_{\Omega_2\times\Omega_2}(f_1-f_2)\,dX\,dY\\
		&\qquad+\left(\frac1{|\Omega_1|}-\frac1{|\Omega_2|}\right)
		\iint_{\Omega_2\times\Omega_2}f_2\,dX\,dY .
	\end{align*}
	The absolute value of the first term is at most $C|(\Omega_1\times\Omega_1)\bigtriangleup(\Omega_2\times\Omega_2)|$. By the bound on the derivative of the integrand with respect to the center,
	\begin{equation*}
		|f_1(X,Y)-f_2(X,Y)|\le C\frac{|R_{\mathrm{center},1}-R_{\mathrm{center},2}|}{R},
	\end{equation*}
	so the second term is at most a constant times this bound. Since $|\Omega_i|\to\pi$, the absolute value of the third term is at most $C\bigl||\Omega_1|-|\Omega_2|\bigr|$. By \eqref{5-1}, the sum of these bounds is at most $Cd$. Thus the right-hand sides of \eqref{5-3} differ by at most $Cb^2d$.

	To compare the left-hand sides of \eqref{5-3}, first subtract the two integrands on $\Omega_1$. For $X\in B_M(0)$, both $r_i(X)/R$ tend uniformly to one. Thus $B''\ge cR^{-1/2}$ on the segment from $r_1(X)$ to $r_2(X)$ by \eqref{4-12}. Hence $B'(r_1(X))-B'(r_2(X))$ has the sign of $R_{\mathrm{center},1}-R_{\mathrm{center},2}$, and
	\begin{equation*}
		\left|\frac1{|\Omega_1|}\int_{\Omega_1}
		\frac{B'(r_1(X))-B'(r_2(X))}{\sqrt R}\,dX\right|
		\ge c\frac{|R_{\mathrm{center},1}-R_{\mathrm{center},2}|}{R}.
	\end{equation*}
	We then compare the averages over $\Omega_1$ and $\Omega_2$ with $r_2(X)$ in both integrands. We may subtract the constant $B'(R_{\mathrm{center},2})/\sqrt R$ from both integrands. For $X\in B_M(0)$,
	\begin{equation*}
		\frac{|B'(r_2(X))-B'(R_{\mathrm{center},2})|}{\sqrt R}
		\le CR^{-1} b|X_1|\le Cb/R
	\end{equation*}
	by the upper bound on $B''$ in \eqref{4-12}. The difference due to integrating over $\Omega_1$ and $\Omega_2$ is bounded by $C(b/R)|\Omega_1\bigtriangleup\Omega_2|$. The difference due to the factors $1/|\Omega_i|$ is bounded by $C(b/R)\bigl||\Omega_1|-|\Omega_2|\bigr|$. Thus the absolute value of the difference between the two left-hand sides of \eqref{5-3} is at least
	\begin{equation*}
		c\frac{|R_{\mathrm{center},1}-R_{\mathrm{center},2}|}{R}
		-C\frac bR|\Omega_1\bigtriangleup\Omega_2| .
	\end{equation*}
	Comparing with $Cb^2d$ on the right and using $b/R=o(b^2)$ and \eqref{5-1} proves \eqref{5-2}.
\end{proof}

Before subtracting the boundary equations, we estimate the difference of the logarithmic potentials of $\Omega_1$ and $\Omega_2$ at the same point on $\partial\Omega_2$.

\begin{lemma}\label{lem5-2}
	For $X=s_2(\theta)e(\theta)$, we have
	\begin{equation}\label{5-4}
		\begin{split}
			&\left\|\int_{\Omega_1}\ln\frac1{|X-Y|}\,dY-\int_{\Omega_2}\ln\frac1{|X-Y|}\,dY
			-\int_{\mathbb T}\ln\frac1{|e(\theta)-e(\eta)|}h(\eta)\,d\eta\right\|_{L^2(\mathbb T)}\\
			&\qquad\le C\delta\|h\|_{L^2(\mathbb T)} .
		\end{split}
	\end{equation}
\end{lemma}

\begin{proof}
In polar coordinates, the difference is the integral between $s_2(\eta)$ and $s_1(\eta)$ for each $\eta$,
	\begin{equation}\label{5-5}
		\begin{split}
			&\int_{\Omega_1}\ln\frac1{|X-Y|}\,dY-\int_{\Omega_2}\ln\frac1{|X-Y|}\,dY\\
			&\quad=\int_0^{2\pi}\int_{s_2(\eta)}^{s_1(\eta)}
			s\ln\frac1{|s_2(\theta)e(\theta)-s e(\eta)|}\,ds\,d\eta .
		\end{split}
	\end{equation}
	In the inner integral set $s=s_2(\eta)+\tau h(\eta)$, $0\le\tau\le1$, so that $ds=h(\eta)\,d\tau$ even when $h(\eta)<0$. For $\theta\neq\eta$,
	\begin{equation*}
		|s_2(\theta)e(\theta)-s e(\eta)|^2
		=(s_2(\theta)-s)^2+s_2(\theta)s|e(\theta)-e(\eta)|^2
	\end{equation*}
	and $|s_2(\theta)-1|+|s-1|\le 2\delta$, so
	\begin{equation*}
		1-C\delta\le
		\frac{|s_2(\theta)e(\theta)-s e(\eta)|^2}{|e(\theta)-e(\eta)|^2}
		\le1+C\delta+\frac{C\delta^2}{|e(\theta)-e(\eta)|^2}.
	\end{equation*}
	Consequently,
	\begin{equation*}
		\left|\ln\frac{|s_2(\theta)e(\theta)-s e(\eta)|}
		{|e(\theta)-e(\eta)|}\right|
		\le C\delta+\frac{1}{2}\ln\left(1+\frac{C\delta^2}{|e(\theta)-e(\eta)|^2}\right).
	\end{equation*}
	Subtracting the kernel on the unit circle from the integrand in \eqref{5-5}, we have
	\begin{align*}
		&s\ln\frac1{|s_2(\theta)e(\theta)-s e(\eta)|}-\ln\frac1{|e(\theta)-e(\eta)|}\\
		&\quad=(s-1)\ln\frac1{|e(\theta)-e(\eta)|}
		+s\ln\frac{|e(\theta)-e(\eta)|}{|s_2(\theta)e(\theta)-s e(\eta)|}.
	\end{align*}
	By the preceding bound on the logarithm and $|s-1|\le\delta$, the absolute value of this difference is bounded by
	\begin{equation*}
		W_\delta(\theta-\eta):=C\delta\left(1+\left|\ln\frac1{|e(\theta)-e(\eta)|}\right|\right)
		+C\ln\left(1+\frac{C\delta^2}{|e(\theta)-e(\eta)|^2}\right).
	\end{equation*}
	This bound depends only on $\theta-\eta$ and is uniform in $\tau$. The first term has $L^1(\mathbb T)$ norm at most $C\delta$. If $t\in[0,\pi]$ is the angle between $e(\theta)$ and $e(\eta)$, then $|e(\theta)-e(\eta)|=2\sin(t/2)\ge2t/\pi$. Hence the integral of the second term is at most
	\begin{equation*}
		C\int_0^\pi \ln(1+C\delta^2/t^2)\,dt
		\le C\delta\int_0^\infty\ln(1+t^{-2})\,dt
		=C\pi\delta .
	\end{equation*}
	Here the last integral equals $\bigl[t\ln(1+t^{-2})\bigr]_0^\infty+2\int_0^\infty(1+t^2)^{-1}\,dt=\pi$ by integration by parts. By \eqref{5-5} and $ds=h(\eta)\,d\tau$, the absolute value of the error in \eqref{5-4} is bounded pointwise by $\int_{\mathbb T}W_\delta(\theta-\eta)|h(\eta)|\,d\eta$. Changing variables $t=\theta-\eta$ and using the triangle inequality for the integral in $L^2(\mathbb T)$, we find
	\begin{align*}
		\left\|\int_{\mathbb T}W_\delta(t)|h(\cdot-t)|\,dt\right\|_{L^2(\mathbb T)}
		&\le\int_{\mathbb T}W_\delta(t)\|h(\cdot-t)\|_{L^2(\mathbb T)}\,dt\\
		&=\|W_\delta\|_{L^1(\mathbb T)}\|h\|_{L^2(\mathbb T)}
		\le C\delta\|h\|_{L^2(\mathbb T)} .
	\end{align*}
\end{proof}

As in \cite{Qin}, we subtract the two formulas in \eqref{4-24} at $X=s_2(\theta)e(\theta)$. Since $\partial_s\tilde F_*=-1/2$ at $s=1$, the first-order term in the difference of the potentials is $h/2$. We keep the difference of the constants for the estimate of the areas below.

\begin{lemma}\label{lem5-3}
	For the two solutions above,
	\begin{equation}\label{5-6}
		\left\|h-\frac1\pi\int_{\mathbb T}\ln\frac1{|e(\theta)-e(\eta)|}h(\eta)\,d\eta
		-2(c(k,v_1)-c(k,v_2))\right\|_{L^2(\mathbb T)}
		\le C(\delta+b^2)d .
	\end{equation}
\end{lemma}

\begin{proof}
	Fix $\theta$ and take $X=s_2(\theta)e(\theta)$. We compare $\tilde F_1(X)$ with $\tilde F_1(s_1(\theta)e(\theta))=0$ along the radial segment. Since $\tilde F_2(X)=0$, integration along this segment implies
	\begin{equation}\label{5-7}
		\begin{split}
			\tilde F_1(X)-\tilde F_2(X) &=-\int_{s_2(\theta)}^{s_1(\theta)} \nabla\tilde F_1(s e(\theta))\cdot e(\theta)\,ds,\\
			\left|\tilde F_1(X)-\tilde F_2(X)-\frac{1}{2}h(\theta)\right| &\le C\delta |h(\theta)|.
		\end{split}
	\end{equation}
	For $s$ between $s_1(\theta)$ and $s_2(\theta)$,
	\begin{equation*}
		\left|\nabla\tilde F_1(se(\theta))\cdot e(\theta)+\frac{1}{2}\right|
		\le\|\tilde F_1-\tilde F_*\|_{C^1(B_2(0))}
		+\left|\partial_s\tilde F_*(se(\theta))+\frac{1}{2}\right|
		\le C\delta ,
	\end{equation*}
	by \eqref{4-27} and $|s-1|\le\delta$.

	Subtracting \eqref{4-24} for the two solutions at $X=s_2(\theta)e(\theta)$, we have
	\begin{align*}
		\tilde F_1(X)-\tilde F_2(X)
		&=c(k,v_1)-c(k,v_2)+\frac1{2\pi}\int_{\R^2}
		\ln\frac1{|X-Y|}\,(\boldsymbol 1_{\Omega_1}-\boldsymbol 1_{\Omega_2})(Y)\,dY\\
		&\quad+E_1(X)-E_2(X).
	\end{align*}
	By \eqref{5-7} and \eqref{5-4}, we may replace the left-hand side by $h/2$ and the integral on the right by the logarithmic integral on the circle. Each error in $L^2(\mathbb T)$ is bounded by $C\delta\|h\|_{L^2(\mathbb T)}$. Since $|X|<2$, the remaining error is bounded by
	\begin{equation*}
		\|(E_1-E_2)(s_2(\cdot)e(\cdot))\|_{L^2(\mathbb T)}
		\le\sqrt{2\pi}\|E_1-E_2\|_{L^\infty(B_2(0))}
		\le C\frac{b\ln R}{R}d,
	\end{equation*}
	where we used \eqref{4-26} and \eqref{5-1}. By the definition of $b$,
	\begin{equation*}
		\frac{b\ln R}{R b^2}=\sqrt{\frac{15}{8}}\,\frac{(\ln R)^{3/2}}R
		\to0.
	\end{equation*}
	Using $\|h\|_{L^2(\mathbb T)}\le d$, we move the logarithmic term and the constant to the left and multiply by $2$ to obtain \eqref{5-6}.
\end{proof}

Let $P$ be the orthogonal projection in $L^2(\mathbb T)$ onto the cosine modes of order $m\ge2$. We estimate these modes directly from the boundary equation.

\begin{lemma}\label{lem5-4}
	For the two solutions above,
	\begin{equation}\label{5-8}
		\|Ph\|_{L^2(\mathbb T)}\le C(\delta+b^2)d .
	\end{equation}
\end{lemma}

\begin{proof}
	For $0<\sigma<1$, the real part of the power series for $-\ln(1-\sigma e^{\mathrm i(\eta-\theta)})$ is
	\begin{equation*}
		\ln\frac1{|1-\sigma e^{\mathrm i(\eta-\theta)}|}
		=\sum_{j=1}^\infty\frac{\sigma^j}{j}\cos(j(\eta-\theta)).
	\end{equation*}
	The mean of each term is zero. For $m\ge1$, we obtain by the cosine addition formula and orthogonality
	\begin{equation*}
		\int_{\mathbb T}\cos(j(\eta-\theta))\cos(m\eta)\,d\eta
		=\begin{cases}
			\pi\cos(m\theta),& j=m,\\
			0,& j\neq m .
		\end{cases}
	\end{equation*}
	Thus the integral of the expansion against $\cos(m\eta)$ is $\pi\sigma^m\cos(m\theta)/m$. To pass to $\sigma=1$, note that for $\sigma\ge1/2$,
	\begin{equation*}
		|1-\sigma e^{\mathrm i(\eta-\theta)}|^2=(1-\sigma)^2+
		\sigma|e(\theta)-e(\eta)|^2\ge\frac{1}{2}|e(\theta)-e(\eta)|^2 .
	\end{equation*}
	By dominated convergence with majorant $C+\left|\ln\frac1{|e(\theta)-e(\eta)|}\right|$, we obtain, for $m\ge1$,
	\begin{equation}\label{5-9}
		\begin{split}
		\int_{\mathbb T}\ln\frac1{|e(\theta)-e(\eta)|}\,d\eta&=0,\\
		\frac1{\pi}\int_{\mathbb T}\ln\frac1{|e(\theta)-e(\eta)|}\cos(m\eta)\,d\eta
		&=\frac1m\cos(m\theta).
		\end{split}
	\end{equation}
	See \cite{Hmi,HM,CJ2} for related Fourier calculations in the study of vortex patches.

	By \eqref{5-9}, if $h=\cos(m\theta)$ for $m\ge1$, the difference between $h$ and the logarithmic integral in \eqref{5-6} equals $(1-1/m)\cos(m\theta)$. Applying $-\partial_\theta$ to $(1-1/m)\cos(m\theta)$, we obtain $(m-1)\sin(m\theta)$, as in the differentiated contour equation in \cite{Qin}. For $m\ge2$ we have $1-1/m\ge1/2$. Since $P$ maps constants to zero, Parseval's identity and \eqref{5-6} imply
	\begin{equation*}
		\frac{1}{2}\|Ph\|_{L^2(\mathbb T)}
		\le\left\|P\left[h-\frac1\pi
		\int_{\mathbb T}\ln\frac1{|e(\theta)-e(\eta)|}h(\eta)\,d\eta\right]\right\|_{L^2(\mathbb T)}
		\le C(\delta+b^2)d .
	\end{equation*}
\end{proof}

To estimate the mean of $h$, we use the relation between the area and the constant term in \eqref{4-24}. At $m=1$, $1-1/m=0$, so we estimate the coefficient of $\cos\theta$ using $\int_{\Omega_i}X_1\,dX=0$ for the two cross-sections.

\begin{lemma}\label{lem5-5}
	For the two solutions above,
	\begin{equation}\label{5-10}
		\left|\int_{\mathbb T}h(\theta)\,d\theta\right|\le C(\delta+b^2)d,\qquad
		\left|\int_{\mathbb T}h(\theta)\cos\theta\,d\theta\right|
		\le C\delta\|h\|_{L^2(\mathbb T)} .
	\end{equation}
\end{lemma}

\begin{proof}
	By symmetry and \eqref{5-9}, the logarithmic kernel has integral zero in $\theta$ for each $\eta$. Integrating the expression inside the norm in \eqref{5-6} with respect to $\theta$ and applying Fubini's theorem and the Cauchy--Schwarz inequality, we have
	\begin{align*}
		\left|\int_{\mathbb T}h(\theta)\,d\theta-4\pi(c(k,v_1)-c(k,v_2))\right|
		&\le C(\delta+b^2)d,\\
		|c(k,v_1)-c(k,v_2)|&\le C\bigl(\|h\|_{L^2(\mathbb T)}+(\delta+b^2)d\bigr)\le Cd.
	\end{align*}
	By the formula for $c(k,v_i)$ in \eqref{4-24},
	\begin{align*}
		c(k,v_1)-c(k,v_2) &=\left(\frac1{2\pi}\ln\frac{R_{\mathrm{center},1}}{b}+c_0\right) (|\Omega_1|-|\Omega_2|)\\
		&\quad+\frac{|\Omega_2|}{2\pi}\ln\frac{R_{\mathrm{center},1}}{R_{\mathrm{center},2}}
		-\frac{k}{15b^2}\bigl(R_{\mathrm{center},1}^{-2}-R_{\mathrm{center},2}^{-2}\bigr).
	\end{align*}
	Since $\ln(R_{\mathrm{center},1}/b)=\ln R+\ln(R_{\mathrm{center},1}/R)-\ln b$ with $\ln(R_{\mathrm{center},1}/R)\to0$ and $-\ln b=O(\ln\ln R)$, the coefficient of $|\Omega_1|-|\Omega_2|$ satisfies, for large $R$,
	\begin{equation*}
		\frac1{2\pi}\ln\frac{R_{\mathrm{center},1}}{b}+c_0
		\ge \frac{1}{4\pi}\ln R .
	\end{equation*}
	For $R/2\le t\le2R$, we have $|(\ln t)'|=1/t\le 2/R$ and $|(kt^{-2})'|=2k/t^3\le C/R$ since $k=3R^2$. The mean value theorem bounds the sum of the absolute values of the last two terms by $Cb^{-2}|R_{\mathrm{center},1}-R_{\mathrm{center},2}|/R$. Dividing by the coefficient of $|\Omega_1|-|\Omega_2|$ and using $(\ln R)^{-1}=(15/8)b^2$, we obtain
	\begin{equation*}
		\bigl||\Omega_1|-|\Omega_2|\bigr|
		\le Cb^2|c(k,v_1)-c(k,v_2)|+C\frac{|R_{\mathrm{center},1}-R_{\mathrm{center},2}|}{R}\le Cb^2d
	\end{equation*}
	by \eqref{5-2}. In polar coordinates, the difference of the areas is
	\begin{align*}
		|\Omega_1|-|\Omega_2|
		&=\frac12\int_{\mathbb T}(s_1(\theta)^2-s_2(\theta)^2)\,d\theta\\
		&=\int_{\mathbb T} h(\theta)\,d\theta
		+\frac12\int_{\mathbb T}(s_1(\theta)+s_2(\theta)-2)h(\theta)\,d\theta .
	\end{align*}
	Since $\|s_i-1\|_{L^\infty(\mathbb T)}<\delta$, we obtain $\left|\int_{\mathbb T}h(\theta)\,d\theta\right|\le C(\delta+b^2)d$.

	To estimate the coefficient of $\cos\theta$, we write the first moment of each core in polar coordinates,
	\begin{equation*}
		0=\int_{\Omega_i}X_1\,dX=\frac13\int_{\mathbb T}s_i(\theta)^3\cos\theta\,d\theta.
	\end{equation*}
	Subtracting the two identities, we have
	\begin{align*}
		0&=\int_{\mathbb T}(s_1^3-s_2^3)\cos\theta\,d\theta\\
		&=3\int_{\mathbb T}h\cos\theta\,d\theta
		+\int_{\mathbb T}h(s_1^2+s_1s_2+s_2^2-3)\cos\theta\,d\theta .
	\end{align*}
	Using $\|s_i-1\|_{L^\infty(\mathbb T)}<\delta$ again, we find $\left|\int_{\mathbb T}h\cos\theta\,d\theta\right|\le C\delta\|h\|_{L^2(\mathbb T)}$.
\end{proof}

\begin{proof}[Proof of Theorem~\ref{thm1-2}]
	The bound on $D$, the limit of $a^2\ln(R_{\mathrm{center}}/a)$ and convergence to the disk follow from Propositions~\ref{pro4-3}, \ref{pro4-4} and \ref{pro4-7}, respectively. Since $D$ lies off the axis and its boundary is a graph in polar coordinates about its center, its revolution is a solid torus.

	For uniqueness, use the evenness of $h$ and the estimates in Lemmas~\ref{lem5-1}, \ref{lem5-4} and \ref{lem5-5} to obtain
	\begin{align*}
		d &\le \frac{|R_{\mathrm{center},1}-R_{\mathrm{center},2}|}{R}
		+\frac{\left|\int_{\mathbb T}h\,d\theta\right|}{\sqrt{2\pi}}
		+\frac{\left|\int_{\mathbb T}h\cos\theta\,d\theta\right|}{\sqrt\pi}
		+\|Ph\|_{L^2(\mathbb T)}\\
		&\le C(\delta+b^2)d .
	\end{align*}
	We choose $\delta$ small and then $K_\infty$ large so that $C(\delta+b^2)<1$ for $k\ge K_\infty$. Thus $d=0$. The centers agree and $h$ is continuous with zero $L^2(\mathbb T)$ norm, so $D_1=D_2$ by \eqref{4-29}. The potential formula \eqref{2-3} then implies $v_1=v_2$ for $r>0$, and continuity extends the equality to the axis.

	For every $k\ge K_\infty$, Corollary~\ref{cor2-8} provides a solution on the Amick--Turner continuum. After translation in $z$, it is the centered solution $v_k$. If $k_n\to k$ in $[K_\infty,\infty)$, then by Lemma~\ref{lem2-6}, $(v_{k_n})$ is relatively compact in $D^{1,2}(\R^5)\cap C^{1,\beta}(\R^5)$ for every $0<\beta<1$. Every convergent subsequence has limit $v_k$ by uniqueness at flux $k$. Thus the entire sequence converges to $v_k$ in this space.
\end{proof}

\subsection*{Data availability}
No data were used for the research described in this article.

\subsection*{Competing interests}
The authors declare that they have no competing interests.

\subsection*{AI assistance statement}
We continue our series of studies on steady vortex rings for the three-dimensional Euler equations. We developed the initial ideas and worked out the first versions of the proofs through discussions with GPT-5.6 Sol.  We subsequently used GPT-6 Astra to simplify parts of the proofs, check references, draw pictures and polish the writing. We independently verified all mathematical statements and proofs and take full responsibility for the content of the manuscript.

\end{document}